\documentstyle[11pt,amsmath,amssymb,epsf,epic]{article}

\numberwithin{equation}{section}

\newcommand{\fa}{{\frak A}}
\newcommand{\fm}{{\frak M}}
\newcommand{\un}{{\mathbb I}}
\newcommand{\ra}{\rightarrow}
\newcommand{\tr}{{\rm Tr}}

\newcommand{\bra}{\langle} 
\newcommand{\ket}{\rangle}

\newcommand{\be}{\begin{equation}}
\newcommand{\ee}{\end{equation}}
\newcommand{\bea}{\begin{eqnarray}}
\newcommand{\eea}{\end{eqnarray}}

\newcommand{\eps}{\epsilon}

\newcommand{\cs}{{\cal S}}

\newcommand{\mqn}{M_{Q_n}}
\newcommand{\e}{{\rm e}}
\newcommand{\p}{{\rm p}}
\renewcommand{\d}{{\rm d}}
\newcommand{\tbeta}{\beta'}
\newcommand{\ext}{{\rm ext}}

\newcommand{\grintl}{[\kern-.18em [}
\newcommand{\grintr}{]\kern-.18em ]}
\newcommand{\ds}{\displaystyle}

\newcounter{resultcounter}[section]

\newtheorem{thm}[resultcounter]{Theorem}
\newtheorem{lem}[resultcounter]{Lemma}
\newtheorem{prop}[resultcounter]{Proposition}

\newtheorem{definition}[resultcounter]{Definition}

\def\bed{\begin{definition}}
\def\eed{\end{definition}}

\def\one{{\mathchoice {\rm 1\mskip-4mu l} {\rm 1\mskip-4mu l} {\rm 1\mskip-4.5mu l} {\rm 1\mskip-5mu l}}}

\def\proof{\noindent{\bf Proof.}\ \ }

 \def\cB{{\cal B}} \def\cC{{\cal C}}
 \def\cE{{\cal E}} \def\cF{{\cal F}}
 \def\cH{{\cal H}} 
  
\def\cM{{\cal M}}  
\def\cP{{\cal P}}  
\def\cS{{\cal S}}

\newcommand{\R}{{\mathbb R}}
\newcommand{\N}{{\mathbb N}}

\newcommand{\C}{{\mathbb C}}
\newcommand{\Z}{{\mathbb Z}}
\newcommand{\E}{{\mathbb E}}
\renewcommand{\P}{{\mathbb P}}

\newcommand{\cme}{{\cal M}_{(E)}}

\def\proof{\noindent{\bf Proof.}\ \ }
\def\qed{\hfill $\Box$\medskip}

\newcommand{\ri}{{\rm i}}
\newcommand{\fer}[1]{(\ref{#1})}
\newcommand{\scalprod}[2]{\left\langle {#1}, {#2}\right\rangle}
\newcommand{\bbbone}{\mathchoice {\rm 1\mskip-4mu l} {\rm 1\mskip-4mu l}
{\rm 1\mskip-4.5mu l} {\rm 1\mskip-5mu l}}
\newcommand{\po}{\overline\omega}

\begin{document}
\title{Infinite Products of Random Matrices\\
and Repeated Interaction Dynamics} 
\author{ Laurent Bruneau\footnote{
Laboratoire AGM, 
Universit\'e de Cergy-Pontoise, Site Saint-Martin, BP 222,
 95302 Cergy-Pontoise, France. Email: laurent.bruneau@u-cergy.fr, http://www.u-cergy.fr/bruneau
}\ , Alain Joye\footnote{
Institut Fourier, UMR 5582, CNRS-Universit\'e de Grenoble I
BP 74, 38402 Saint-Martin d'H\`eres, France. Email: Alain.Joye@ujf-grenoble.fr
}\,\,\footnote{Laboratoire de Physique et Mod\'elisation des Milieux
  Condens\'es, UMR 5493, CNRS-Universit\'e de Grenoble I, BP 166,
  38042 Grenoble, France}\ , Marco Merkli\footnote{Memorial University,
Canada. Email: merkli@math.mun.ca, http://www.math.mun.ca/\ $\widetilde{}$\,merkli/}
}

\date{
}

\maketitle
\vspace{-1cm}
\begin{abstract}
Let $\Psi_n$ be a product of $n$ independent, identically distributed random matrices $M$, with the properties that $\Psi_n$ is bounded in $n$, and that $M$ has a deterministic (constant) invariant vector. Assume that the probability of $M$ having only the simple eigenvalue 1 on the unit circle does not vanish.
 We show that $\Psi_n$ is the sum of a fluctuating and a decaying process. 
The latter converges to zero almost surely, exponentially fast as $n\rightarrow\infty$. The fluctuating part converges in Cesaro mean to a limit that is characterized explicitly by the deterministic invariant vector and the spectral data of ${\mathbb E}[M]$ associated to 1. No additional assumptions are made on the matrices $M$; they may have complex entries and not be invertible. 

We apply our general results to two classes of dynamical systems: inhomogeneous Markov chains with random transition matrices (stochastic matrices), and random repeated interaction quantum systems. In both cases, we prove ergodic theorems for the dynamics, and we obtain the the limit states.
\end{abstract}

\thispagestyle{empty}
\setcounter{page}{1}
\setcounter{section}{1}


\setcounter{section}{0}

\section{Introduction}
\label{sec:intro}

In this paper we study products of infinitely many independent, identically distributed random matrices. The matrices we consider satisfy two basic properties, reflecting the fact that they describe the dynamics of random quantum or classical dynamical systems. The first property is that the norm of any product of such matrices is bounded uniformly in the number of factors. It reflects the fact that the underlying dynamics is in a certain sense norm-preserving. The second property is that there is a deterministic invariant vector. This represents a normalization of the dynamics. Two important examples of systems falling into this category are inhomogeneous Markov chains with {\it random transition matrices}, i.e. products of {\it random stochastic matrices}, as well as {\it repeated interaction open quantum systems}. In this paper we present general results on infinite products of random matrices, and applications of these results to the two classes of dynamical systems mentioned. Our main results are convergence theorems of the infinite random matrix product,  Theorems \ref{thm1}, \ref{thm2} and \ref{thm3}. They translate into ergodic theorems for the corresponding dynamical systems, Theorems \ref{thmmarkov} and \ref{thm4}.


\subsection{General Results}

Let $M(\omega)$ be a random matrix on ${\mathbb C}^d$, with
probability space $(\Omega,\cF,\p)$. We say that $M(\omega)$ is a {\it random reduced dynamics operator} (RRDO) if
\begin{itemize}
\item[(1)] There exists a norm $||| \cdot |||$ on $\C^d$ such that, for all $\omega$,
      $M(\omega)$ is a contraction on $\C^d$ endowed with the norm $|||\cdot|||$. 
\item[(2)] There is a vector $\psi_\cS$, constant in $\omega$, such that $M(\omega)\psi_\cS=\psi_\cS$, for all $\omega$.
\end{itemize}
We shall normalize $\psi_\cS$ such that $\|\psi_\cS\|=1$ where $\|\cdot\|$ denotes the euclidean norm.

\noindent
To an RRDO $M(\omega)$, we associate the (iid) {\it random reduced dynamics process} (RRDP)
$$
\Psi_n(\po):=M(\omega_1)\cdots M(\omega_n),\qquad \po\in{\Omega}^{{\mathbb N}^*}.
$$ 

We show that $\Psi_n$ has a decomposition into an exponentially decaying part and a fluctuating part. To identify these parts, we proceed as follows. It follows from (1) and (2) that the spectrum of an RRDO $M(\omega)$ must lie inside the closed complex unit disk, and that $1$ is an eigenvalue (with eigenvector $\psi_\cS$). Let $P_{1}(\omega)$ denote the spectral projection of $M(\omega)$ corresponding to the eigenvalue $1$ ($\dim P_1(\omega)\geq 1$), and let $P^*_1(\omega)$ be its adjoint operator. Define
\begin{equation}
\psi(\omega) := P_1(\omega)^*\psi_\cS,
\label{m13}
\end{equation}
and set
$$
P(\omega)=|\psi_\cs\rangle \langle\psi(\omega)|.
$$
{}For $\psi,\phi\in{\mathbb C}^d$, we denote by $|\psi\rangle\langle\phi|$ the rank-one operator ${\mathbb C}^d \ni\chi\mapsto |\psi\rangle\langle\phi|\chi = \scalprod{\phi}{\chi}\psi$, and our convention is to take the inner products linear in the second factor. We put 
$$
Q(\omega)=\one-P(\omega).
$$
Note that the vector $\psi(\omega)$ is normalized as $\scalprod{\psi_\cS}{\psi(\omega)}=1$. We decompose $M(\omega)$ as
\begin{equation}
M(\omega) = P(\omega)+Q(\omega)M(\omega)Q(\omega)=: P(\omega) +M_Q(\omega).
\label{m14}
\end{equation}
%
%
Taking into account this decomposition, we obtain (c.f. Proposition \ref{prop:psiform})
\begin{equation}
\Psi_n(\po):=M(\omega_1)\cdots M(\omega_n) = |\psi_\cS\rangle\langle\theta_n(\po)| + M_Q(\omega_1)\cdots M_Q(\omega_n),
\label{m15}
\end{equation}
where $\theta_n(\po)$ is the Markov process
\begin{equation}
\theta_n(\po) =M^*(\omega_n)\cdots M^*(\omega_2)\psi(\omega_1),
\label{m16}
\end{equation}
$M^*(\omega_j)$ being the adjoint operator of $M(\omega_j)$. We analyze the two parts in the r.h.s. of \fer{m15} separately. 

\medskip
Let $\cme$ be the set of RRDOs $M$ whose spectrum on the complex unit circle consists only of a simple eigenvalue $\{1\}$.

On $\Omega^{\N^*}$ we define the probability measure $\d\P$ in a standard fashion by 
$$
\d\P =\Pi_{j\geq 1}\d \p_j, \ \ \ \mbox{where } 
\ \ \ \d \p_j\equiv \d \p, \ \  \forall j\in \N^*.
$$


\begin{thm}[Decaying process]
\label{thm1}
Let $M(\omega)$ be a random reduced dynamics operator. Suppose that $\p(M(\omega)\in\cme)>0$.
Then there exist a set $\Omega_1\subset\Omega^{\N^*}$ and constants $C,\alpha>0$ s.t. ${\mathbb P}(\Omega_1)=1$ and s.t. for any $\po\in\Omega_1$, there exists  an $n_0(\po)$ s.t. for any $n\geq n_0(\po)$,
\begin{equation}
\| M_{Q}(\omega_1)\cdots M_{Q}(\omega_n)\|\leq C\e^{-\alpha n}.
\label{m17}
\end{equation}
\end{thm}

\noindent {\it Remarks.}  1. In the case where $M(\omega)=M$ is constant, and $M\in\cme$, one readily shows that for any $\epsilon>0$ there is a $C_\epsilon$ such that $\|(M_Q)^n\|\leq C_\epsilon\e^{-n(\gamma-\epsilon)}$, for all $n\geq 0$, and where $\gamma=\min_{z\in{\rm spec}(M)\backslash\{1\}}|\, \log|z|\,|$ (see e.g. \cite{bjm}, Proposition 2.2). It is remarkable that in the random case, the mere condition of $M$ having an arbitrarily small, non-vanishing probability to be in $\cme$ suffices to guarantee the exponential decay of the product in \fer{m17}.

2. Any stochastic matrix whose entries are all nonzero belongs to $\cme$.

3. If $\{1\}$ is a simple eigenvalue of $M(\omega)$ then the decomposition
\fer{m14} is just the spectral decomposition of the matrix $M(\omega)$. 

4. The choice \fer{m13} ensures that $\psi(\omega)$ is an eigenvector of $M^*(\omega)$. Other choices of measurable $\psi(\omega)$ which are bounded in $\omega$ lead to different decompositions of $M(\omega)$, and can be useful as well. For instance, if $M(\omega)$ is a bistochastic matrix, then one can take for $\psi(\omega)$ an $M^*(\omega)$-invariant vector which is constant in $\omega$.

\medskip
\noindent
Our next result concerns the asymptotics of the Markov process \fer{m16}. Set ${\mathbb E}[f]=\int_\Omega f(\omega)\d\p(\omega)$ for a random variable $f$, and denote by $P_{1,{\mathbb E}[M]}$ the spectral projection of ${\mathbb E}[M]$ onto the eigenvalue $\{1\}$. 

\begin{thm}[Fluctuating process]
\label{thm2}
Let $M(\omega)$ be a random reduced dynamics operator. Suppose that $\p(M(\omega)\in\cme)>0$. Then we have ${\mathbb E}[M]\in\cme$. Moreover, there exists a set $\Omega_2\subset\Omega^{\N^*}$ s.t. ${\mathbb P}(\Omega_2)=1$ and, for all $\po\in\Omega_2$,
\begin{equation}
\lim_{N\rightarrow\infty}\frac 1N\sum_{n=1}^N\theta_n(\po) = \theta,
\label{m21}
\end{equation}
where 
\begin{equation}
\theta = \left( \bbbone - {\mathbb E}[M_Q]^*\right)^{-1}{\mathbb E}[\psi]=P^*_{1,{\mathbb E}[M]}{\mathbb E}[\psi]=P^*_{1,{\mathbb E}[M]}\psi_\cS.
\label{m21.01}
\end{equation}
\end{thm}

\medskip
\noindent {\it Remarks.} 5. In the case where $M$ is constant in $\omega$, we have ${\mathbb E}[M_Q]^*=(M_Q)^*$, ${\mathbb E}[\psi]=\psi$, and under the assumption of Theorem \ref{thm2}, that $M\in\cme$. Therefore, $P_1=P=|\psi_\cs\rangle\langle\psi|$ and hence $Q^*\psi=0$, and $(M_Q)^*\psi=0$. Consequently, we have $\theta=\psi$. This coincides with the results of \cite{bjm}. As we will see in Section \ref{ssec:ergocvg}, the latter equality is not satisfied for general, $\omega$-dependent matrices $M$.

6. The ergodic average limit of $\theta_n(\po)$ does not depend on the particular choice of $\psi(\omega)$. This follows from the last equality in \fer{m21.01}.

7. We show in Proposition \ref{limpsi} that for every {\it fixed}
$\po$, $\theta_n(\po)$ converges if and only if
$\psi(\omega_n)$ converges, and that the limits coincide if they
exist.

\medskip
\noindent
Combining Theorems \ref{thm1} and \ref{thm2} we obtain the following result.
\begin{thm}[Ergodic theorem for RRDP]
\label{thm3}
Let $M(\omega)$ be a random reduced dynamics operator. Suppose  $\p(M(\omega)\in\cme)>0$. Then there exists a set $\Omega_3\subset\Omega^{\N^*}$ s.t. ${\mathbb P}(\Omega_3)=1$ and, for all $\po\in\Omega_3$,
\begin{equation}
\lim_{N\rightarrow\infty}\frac 1N\sum_{n=1}^N M(\omega_1)\cdots
M(\omega_n) =|\psi_\cS\ket\bra \theta|=P_{1,{\mathbb E}[M]}.
\label{m21.1}
\end{equation}
\end{thm}

\noindent {\it Remark.\ } 8. If one can choose $\psi(\omega)\equiv\psi$ to be independent of $\omega$ (see also remark 4 above), then we show (see (\ref{eq:theta2})) that $\theta_n(\po)=\psi$, for all $n,\po$. It thus follows from \fer{m15}-\fer{m17} that $
\lim_{n\rightarrow\infty} M(\omega_1)\cdots
M(\omega_n)= |\psi_\cS\ket\bra \psi|$, a.s., exponentially fast.


\subsection{Inhomogeneous random Markov chains}

Products of random matrices satisfying conditions (1) and (2) are important in the study of inhomogenous random Markov chains.
A finite Markov chain is a system consisting of $d< \infty$ states on which a
stochastic process is defined in the following way. At time step $n$, 
$n\in{\mathbb N}^*$, the probability to make a transition from state $i$ to 
state $j$, where $i, j\in \{1,\ldots, d\}$, is given by $p_{ij}^{(n)}$. Hence the relation
$\sum_{j=1}^d p_{ij}^{(n)}=1$. Arranged into a $d\times d$ matrix 
$(M_n)_{ij}=p^{(n)}_{ij}$, these probabilities yield a stochastic matrix. The probability of jumping from state $i$ to state $j$ between times 
$0$ and $n$ is given by the matrix element $ij$ of the product
\be\label{imc}
\Psi_n=M_1M_{2}\cdots M_n.
\ee
The Markov chain \fer{imc} is called {\it inhomogeneous} if the transition probabilities depend on the time step $n$. We have that for all $n$,
\begin{itemize}
\item[-] $M_n$ is a contraction for the norm $|||M|||=\max_{i=1,\ldots,d}\sum_{j=1}^d|(M)_{ij}|$, and 

\item[-] the vector $\psi_{\cal S}=\frac{1}{\sqrt{d}}(1,1, \ldots, 1)^T$ is a normalized eigenvector of 
$M_n$, associated to the eigenvalue 1.
\end{itemize}
Consequently, if the transition probabilities are chosen randomly at each time step, the corresponding operator $M(\omega)$ is an RRDO, and the product (\ref{imc}) is an RRDP. If $p^{(n)}_{ij}>0$, for all $i,j$, then the Perron-Frobenius Theorem asserts that 1 is a simple eigenvalue of $M_n$, and that all other eigenvalues lie inside the open complex unit disk. Thus the conclusions of Theorems \ref{thm1} -\ref{thm3} hold for $M(\omega)$ a random stochastic matrix, which, with positive probability, has strictly positive entries. 

Results on convergence of products of random stochastic matrices, or non-negative matrices, 
are numerous, see e.g. the references in Section \ref{relatedworkssection}. However, these results  
mostly concern properties of the limiting distribution, if it exists, in terms 
of the properties of the distribution of the RRDOs. Those results 
rely heavily on the positivity of the elements of the considered matrices. When studying convergence in distribution, the order of the factors in the product
does not matter, and usually products of the form
\be\label{wrong}
\Phi_n=M_nM_{n-1}\cdots M_1
\ee
as studied as well (compare with \fer{imc}). While products of the form \fer{wrong} are not our main concern (and are easier to study), our techniques still yield results for them that are stronger than those for products \fer{imc}. We show  in Section \ref{ssec:met} a general result on products \fer{wrong} which, applied to products of stochastic matrices, yields the following result.


\begin{thm}
\label{thmmarkov}
Suppose that $M(\omega)$ takes values in the set of stochastic matrices,  and that $p(M(\omega)\in {\cal P}_E)>0$, where ${\cal P}_E$
denotes the set of stochastic matrices with stricly positive entries. Then there exists an $\alpha>0$ and a set $\Omega_4\subset{\Omega}^{{\mathbb N}^*}$, such that ${\mathbb P}(\Omega_4)=1$, and, for all $\po\in\Omega_4$,
\be
\Phi_n(\po)=M(\omega_n)M(\omega_{n-1})\cdots M(\omega_1)=
|\psi_S\ket\bra \eta_\infty(\omega)|+O_{\po}(e^{-\alpha n}),
\ee
where $\eta_\infty(\omega)$ is a vector valued random variable, where 
$\psi_S=\frac{1}{\sqrt{d}}(1,1,\cdots, 1)^T$, and where the remainder term depends on 
$\po\in\Omega_4$.
\end{thm}
\noindent
{\it Remark.\ } 9. As expected, ${\mathbb E}[\eta_\infty]=\theta$, with $\theta$ given in \fer{m21.01}.


\subsection{Random repeated interaction open quantum systems}
 
In this section we present the mathematical framework of open quantum systems. We establish the link between the dynamics and infinite products of random matrices, and we present an application of Theorems \ref{thm1}-\ref{thm3}.


\subsubsection{Definition of repeated interaction quantum systems} 

A repeated interaction quantum system consists of a subsystem $\cS$ which interacts successively with the elements $\cE_m$ of a chain  $\cC=\cE_1+\cE_2+\cdots$ of independent quantum systems. At each moment in time, $\cS$ interacts precisely with one $\cE_m$ ($m$ increases as time does), while the other elements in the chain evolve freely according to their intrinsic (uncoupled) dynamics. The complete evolution is described by the intrinsic dynamics of $\cS$ and of all the $\cE_m$, plus an interaction between $\cS$ and $\cE_m$, for each $m$. The latter consists of an interaction time $\tau_m>0$, and an interaction operator $V_m$ (acting on $\cS$ and $\cE_m$); during the time interval $[\tau_1+\cdots+\tau_{m-1}, \tau_1+\cdots+\tau_{m})$, $\cS$ is coupled to $\cE_m$ via $V_m$. 
Systems with this structure are important not only from a mathematical, but also from a physical point of view. They arise in fundamental experiments on the interaction of matter with quantized radiation. We refer to the end of this introduction for more information and references concerning this aspect of our work.

According to the fundamental principles of quantum mechanics, states
of the systems $\cS$ and $\cE_m$ are given by normalized vectors (or
density matrices) on Hilbert spaces $\cH_\cS$ and $\cH_{\cE_m}$,
respectively, \cite{AJP,BR}\footnote{A normalized vector $\psi$ defines a ``pure'' state $A\mapsto \scalprod{\psi}{A\psi}={\rm Tr}(\varrho_\psi A)$, where $\varrho_\psi=|\psi\rangle\langle\psi|$. A general ``mixed'' state is given by a density matrix $\varrho =\sum_{n\geq 1}p_n\varrho_{\psi_n}$, where the probabilities $p_n\geq 0$ sum up to one, and where the $\psi_n$ are normalized vectors.}. We assume that $\dim\cH_\cS<\infty$, while $\dim\cH_{\cE_m}$ may be infinite. Observables $A_\cS$ and $A_{\cE_m}$ of the systems $\cS$ and $\cE_m$ are bounded operators forming {\it
  von Neumann algebras} $\fm_\cS\subset \cB(\cH_\cS)$ and
$\fm_{\cE_m}\subset \cB(\cH_{\cE_m})$. They evolve according to the {\it Heisenberg
  dynamics} ${\mathbb R}\ni t\mapsto \alpha^t_\cS(A_\cS)$ and ${\mathbb
  R}\ni t\mapsto \alpha^t_{\cE_m}(A_{\cE_m})$, where
$\alpha^t_\cS$ and $\alpha^t_{\cE_m}$ are $*$-automorphism groups of
$\fm_\cS$ and $\fm_{\cE_m}$, respectively, see
e.g. \cite{BR}. We introduce distinguished {\it reference states}, given by vectors $\psi_{\cS}\in\cH_\cS$ and $\psi_{\cE_m}\in \cH_{\cE_m}$. Typical choices for $\psi_\cS$, $\psi_{\cE_m}$ are equilibrium (KMS) states for the dynamics $\alpha^t_\cS$, $\alpha^t_{\cE_m}$, at inverse temperatures $\beta_\cS$, $\beta_{\cE_m}$. The Hilbert space of states of the total system is the tensor product 
$$
\cH=\cH_\cS\otimes\cH_\cC,
$$
where $\cH_\cC=\bigotimes_{m\geq 1}\cH_{\cE_m}$, and where the infinite product is taken with respect to $\psi_\cC=\bigotimes_{m\geq 1}\psi_{\cE_m}$. The non-interacting dynamics is the product of the individual dynamics, defined on the algebra $\fm_{\cS}\bigotimes_{m\geq 1}\fm_{\cE_m}$ by $\alpha_\cs^t\bigotimes_{m\geq 1}\alpha_{\cE_m}^t$. It is useful to consider the dynamics in the {\it Schr\"odinger picture}, i.e. as acting on vectors in $\cH$. To do this, we first implement the dynamics via unitaries, satisfying 
\begin{equation}
\alpha_\#^t(A_\#) = \e^{\ri t L_\#} A_\# \e^{-\ri t L_\#},\ t\in{\mathbb R},\ \ \mbox{and $L_\#\psi_\#=0$},
\label{m1}
\end{equation}
for any $A_\#\in {\frak M}_\#$, where $\#$ stands for either $\cS$ or $\cE_m$. The self-adjoint operators $L_\cS$ and $L_{\cE_m}$, called {\it Liouville operators}, act on $\cH_\cS$ and $\cH_{\cE_m}$, respectively. The existence and uniqueness of $L_\#$ satisfying \fer{m1} is well known, under general assumptions on the reference states $\psi_\#$ \cite{BR}. In particular, \fer{m1} holds if the reference states are equilibrium states. Let $\tau_m>0$ and $V_m\in\fm_{\cS}\otimes\fm_{\cE_m}$ be the interaction time and interaction operator associated to $\cS$ and $\cE_m$. We define the (discrete) repeated interaction Schr\"odinger dynamics of a state vector $\phi\in\cH$, for $m\geq 0$, by 
\begin{equation}
U(m)\phi = \e^{-\ri \widetilde L_m}\cdots\e^{-\ri \widetilde L_2}\e^{-\ri \widetilde L_1}\phi,
\label{m2}
\end{equation}
where 
\begin{equation}
\widetilde L_k = \tau_k L_k+\tau_k \sum_{n\neq k} L_{\cE_n}
\label{m3}
\end{equation}
describes the dynamics of the system during the time interval $[\tau_1+\cdots+\tau_{k-1},\tau_1+\cdots+\tau_k)$, which corresponds to the time-step $k$ of our discrete process. Hence $L_k$ is
\begin{equation}
L_k=L_\cS +L_{\cE_k} + V_k,
\label{m4}
\end{equation}
acting on $\cH_\cS\otimes\cH_{\cE_k}$. Of course, we understand that the operator $L_{\cE_n}$ in \fer{m3} acts nontrivially only on the $n$-th factor of the Hilbert space $\cH_\cC$ of the chain. 

A state $\varrho(\cdot)={\rm Tr}(\rho\, \cdot\,)$ given by density matrix $\rho$ on $\cH$ is called a {\it normal state}. (A density matrix is a self-adjoint, non-negative trace-class operator of unit trace.) Our goal is to understand the large-time asymptotics ($m\rightarrow \infty$) of expectations
\begin{equation}
\varrho\left(U(m)^* O U(m)\right)=\varrho(\alpha_{RI}^m(O)),
\label{m5}
\end{equation}
for normal states $\varrho$ and certain classes of observables $O$ that we specify below. We denote the random repeated interaction dynamics by
\begin{equation}
\alpha_{RI}^m(O) = U(m)^* O U(m).
\label{m5.1}
\end{equation}


\subsubsection{Reduced dynamics and random matrix products} 

Let us explain how we link the dynamics to a product of reduced dynamics operators. In order not to muddle the essence of our argument, we only consider the expectation of an observable $A_\cS\in\fm_\cS$, and we take the initial state of the entire system to be given by the vector \begin{equation}
\psi_0=\psi_\cS\otimes\psi_\cC, 
\label{m0}
\end{equation}
where the $\psi_\cS$ and $\psi_\cC$ are the reference states introduced above. (See Section \ref{sec:reduction} for the general case.) The expectation of $A_\cS$ at the time-step $m$ is 
\begin{equation}
\scalprod{\psi_0}{\alpha^m_{RI}(A_\cS)\psi_0} = \scalprod{\psi_0}{P\e^{\ri\widetilde L_1}\cdots \e^{\ri\widetilde L_m}A_\cS\, \e^{-\ri\widetilde L_m}\cdots \e^{-\ri\widetilde L_1}P\psi_0}.
\label{m6}
\end{equation}
We write simply $A_\cS$ for $A_\cS\otimes\bbbone_\cC$, and we have introduced  
\begin{equation}
P=\bbbone_{\cH_\cS}\bigotimes_{m\geq 1} P_{\psi_{\cE_m}},
\label{m7}
\end{equation}
the orthogonal projection onto $\cH_\cS\otimes {\mathbb C} \psi_\cC$. A {\it first important ingredient} of our approach is to construct operators $K_k$ with the properties
\begin{eqnarray}
\e^{\ri\widetilde L_k} A \e^{-\ri\widetilde L_k} &=&\e^{\ri K_k} A \e^{-\ri K_k},\label{m8}\\
K_k\,\psi_\cS\otimes \psi_\cC&=&0, \label{m8.1}
\end{eqnarray}
where $A$ in \fer{m8} is any observable of the total system. Equality \fer{m8} means that the operators $K_k$ implement the same dynamics as the $\widetilde L_k$. Relation \fer{m8.1} selects a unique generator of the dynamics among all that satisfy \fer{m8}. The existence of operators $K_k$ satisfying \fer{m8} and \fer{m8.1} is linked to the deep Tomita-Takesaki theory of von Neumann algebras, c.f. \cite{bjm} and references therein. It turns out that the $K_k$ are non-normal operators on $\cH$ (while the $\widetilde L_k$ are self-adjoint). We combine \fer{m8} with \fer{m6} to obtain 
\begin{equation}
\scalprod{\psi_0}{\alpha_{RI}^m(A_\cS)\psi_0} = \scalprod{\psi_0}{P\e^{\ri K_1}\cdots\e^{\ri K_m} P A_\cS\,\psi_0}.
\label{m9}
\end{equation}
A {\it second important ingredient} of our approach is to realize that
\begin{equation}
P\e^{\ri K_1}\cdots\e^{\ri K_m} P = P\e^{\ri K_1}P \cdots P\e^{\ri K_m}P,
\label{m10}
\end{equation}
which follows from the independence of the systems $\cE_m$ (see \fer{mm100}).  We identify $P\e^{\ri K_k}P$ with an operator $M_k$ on $\cH_\cS$, and thus obtain, from \fer{m9} and \fer{m10},
\begin{equation}
\scalprod{\psi_0}{\alpha_{RI}^m(A_\cS)\psi_0} = \scalprod{\psi_\cS}{M_1\cdots M_m A_\cS\, \psi_\cS}.
\label{m11}
\end{equation}
Because the operators $M_k=P\e^{\ri K_k}P$ implement a dynamics, we
can show (Lemma \ref{contraction}) that the $M_k$ are contractions for some
suitable norm $|||\cdot |||$ of ${\mathbb C}^d$. Moreover, it follows from \fer{m8.1}
that $M_k\psi_\cS=\psi_\cS$, for all $k$. The matrices $M_k$ satisfy conditions (1) and (2) above and are thus RRDOs.

\medskip
\noindent
{\it Remark.\ } 10. The process of reduction of the dynamics explained in this section {\it does not involve any approximation}. This is in contrast e.g. with master equation techniques, where a reduced dynamics is obtained using Born and Markov approximations, or so-called weak coupling approximations. Our method allows us to use the structure of repeated interaction systems, and to do without all these approximations.


\subsubsection{Ergodicity of repeated interaction systems}

Let us denote by $\alpha^{n,\po}_{RI}$, $\po\in\Omega^{{\mathbb N}^*}$, the random reduced operator process obtained from \fer{m5.1}, \fer{m11}, where the $M_j$ in \fer{m11} are iid random matrices. We call $\alpha_{RI}^{n,\po}$ the random repeated interaction dynamics determined by the RRDO $M(\omega)=P\e^{\ri K(\omega)}P$.

\begin{thm}[Ergodic theorem for random repeated interaction systems]\ 
\label{thm4}
{}\ Let $\alpha_{RI}^{n,\po}$ be the random repeated interaction dynamics determined by an RRDO $M(\omega)$. Suppose that $\p(M(\omega)\in\cme)>0$. Then there exists a set $\Omega_5\subset \Omega^{\N^*}$, s.t.  ${\mathbb P}(\Omega_5)=1$, and s.t. for any $\po\in\Omega_5$, any normal state $\varrho$ and any $A_\cs\in\fm_\cS$,
\begin{equation}
\lim_{N\rightarrow\infty} \frac 1N\sum_{n=1}^N \varrho\left( \alpha^{n,\po}_{RI}(A_\cS)\right)=\scalprod{\theta}{A_\cS \psi_\cS},
\label{m22}
\end{equation}
where $\theta$ is given by \fer{m21.01}.
\end{thm}

\medskip
\noindent
{\it Remarks.\ }  11. For the normal state $\varrho$ given by the density matrix $\rho=|\psi_0\rangle\langle \psi_0|$, \fer{m22} follows readily from Theorem \ref{thm3} and relation \fer{m11}. To prove the result for arbitrary normal states, one uses a property of the reference state $\psi_0$ called {\it cyclicity and separability}, \cite{BR,bjm}.

12. Our setup allows us to treat systems having various sources of randomness. Not only do we allow for random interactions, but also for random characteristics of the systems $\cE_m$ and $\cS$ (interesting cases are random temperatures and dimensions of the $\cE_m$ and of $\cS$).

13. In \cite{bjm2}, we prove a similar result for a general class of so-called {\it instantaneous observables} (involving also operators on the chain $\cC$), and we discuss the physical properties of the asymptotic state. See also \cite{bjm} for the case where $M(\omega)\equiv M$.

14. As mentioned above, $\scalprod{\psi_0}{\alpha^{n,\po}_{RI}(A_\cS)\psi_0}$ is fluctuating, it does not converge pointwise (in $\po \in\Omega^{{\mathbb N}^*}$), but only in the ergodic average sense.

15. We present in Section \ref{sec:spin} the explicit example of a spin-spin repeated interaction open quantum system.


\subsection{Related works}
\label{relatedworkssection}

\qquad 
{\it Random matrices.\ }
Results on the convergence (of some kind) of 
products of random stochastic matrices, and non-negative matrices, 
are numerous, see e.g. \cite{ks,m,s} and references therein. These results  
mostly concern properties of the limiting distribution, if it exists, in terms 
of the properties of the distribution of the RRDOs. The techniques used to obtain those results 
rely heavily on the positivity of matrix elements. Random matrix products have been heavily studied also for matrices in $Gl_d({\mathbb R})$, see e.g.  
\cite{g,k,la}. Again, the main focus of these works is on the study
of the properties of the limiting distribution, if it exists, and on the 
properties of the Lyapunov exponents. In this case, the group property
of the set of invertible matrices  $Gl_d({\mathbb R})$ plays a prominent role in the derivation of the results. 

By contrast, besides conditions (1) and (2) defining
RRDO, we do not require our matrices to be real valued, to have positive entries, or to be invertible. 
Moreover, we are  concerned here with the limiting properties of the products only, not
with the limiting distribution. On the one hand, we get a.s. convergence 
results for the products $\Phi_n$ (Theorem \ref{thmmarkov}).  On the other hand,
in order to eliminate the unavoidable fluctuations in the products $\Psi_n$, 
we resort to a limit in an average sense, the 
Cesaro limit (Theorems \ref{thm2}-\ref{thm3}). Let us also point out that our results show as well that the top Lyapunov exponent of products of RRDOs is zero, and that it is almost surely of multiplicity one, see Theorem \ref{prop:randlyap}.


Getting informations on the fluctuations of the process around its limiting value is certainly an interesting and important issue. It amounts to getting informations about the law of the vector valued random variable $\eta_\infty$ of Theorem \ref{thmmarkov}, which is quite difficult in general. There are recent partial results about aspects of the law of such random vectors in case they are obtained by means of matrices belonging to some subgroups of $Gl_d(\mathbb R)$ satisfying certain irreducibility conditions, see e.g. \cite{dsetal}. However, these results do not apply to our situation.

Random ergodic theorems for products $\Psi_n(\po)$ have been obtained in a more general framework in \cite{beckschwarz}. They prove the almost sure existence of a Cesaro limit for these products. However, the identification of the limit is not provided by this general result. This identification relies on the detailed properties of the matrices involved, and in particular in the separation of a fluctuating part from a decaying part in the dynamical process. In this respect, our contribution consists in identifying completely the Cesaro limit of products of RRDO.


Actually, in the case under consideration, Theorem \ref{thm1} yields enough information so that Theorem \ref{thm3} becomes a consequence of it and of the existence result proven in \cite{beckschwarz}. Since a direct self-contained proof of Theorem \ref{thm3} is easy to derive, given Theorem \ref{thm1}, we present it for completness.

We close this comparison by noting that the paper \cite{ben}
contains results, in a deterministic setup, which are close to some of the results
we prove in Section 2. The authors of \cite{ben} consider infinite products of  
contractions for some norm on $M_d({\mathbb C})$, in the reverse order with respect 
to RRDP. They prove convergence of the product under the following assumptions:
i) there exists a subsequence of the matrices appearing in the product which
converges to a matrix which is paracontracting; ii) the set of invariant vectors
of this matrix is contained in the set of invariant vectors of all matrices 
appearing in the product. This last assumption implies
our condition (2) in a deterministic setup. However, the products 
in \cite{ben} are in the reverse order with respect to those we address, and they address the deterministic framework only.

{\it Repeated interaction systems.\ } In the experimental setup of a ``One-Atom Maser'' the system $\cS$ represents one or several modes of the quantized electromagnetic field in a cavity, while the elements $\cE_m$ describe atoms injected into the cavity, one by one, interacting with the radiation while passing through the cavity. After interaction, the atoms encode certain properties of the field that can be measured after they exit the cavity, \cite{MWM,WVHW}. In an idealized model for this process, one assumes that all the elements $\cE_m$ represent a copy of the same, fixed quantum system, and that the interaction is given by a fixed interaction time and a fixed interaction operator (the same for all $m$). Such {\it idealized repeated interaction systems} have been analyzed mathematically in \cite{WBKM,bjm}. However, it is clear that in actual experiments, neither the interaction time $\tau_m$ (or $V_m$) nor the elements $\cE_m$ can be exactly the same for all $m$! Typically, the interaction time will be {\it random}, given e.g. by a Gaussian or by a uniform distribution around a mean value, and the state of the incoming atoms will be random as well, for instance determined by a temperature that fluctuates slightly around a mean temperature. (In experiments, the atoms are ejected from an atom oven, then they are cooled down to a wanted temperature before entering the cavity.) It is therefore important to develop a theory that allows for {\it random repeated interactions}, as we do in the present paper. Our approach is based on the techniques of non-equilibrium quantum statistical mechanics developed in \cite{bjm}. We are not aware of any other work dealing with variable and random interactions in the quantum context.


\subsection{Organization of this paper}

Our paper is organized as follows. In Section \ref{sec:determinist},
we prove several deterministic results, including the fact that $\Psi_n$ converges (for fixed $\omega$) if and only if there are no fluctuations in the matrix product (Proposition \ref{limpsi}). In Section \ref{sec:random} we prove our main results, Theorems \ref{thm1} and \ref{thm2}.  In Section \ref{sec:reduction} we explain the link between repeated interaction quantum dynamics and matrix products. We apply, in Section \ref{sec:spin}, our general results to an explicit quantum dynamical system, a random repeated interaction spin system.


\section{Deterministic results}
\label{sec:determinist}

In this section, we derive some algebraic formulae and some uniform
bounds which will play a crucial role in our analysis. However, there
is no probabilistic issue involved here and we shall therefore {\it
  freeze} the random variable. We will thus simply denote
$M_j=M(\omega_j)$, and   
\be\label{def:psi}
\Psi_n:= M_1 \cdots M_n.
\ee


\subsection{Decomposition of the $M_j$}
\label{ssec:notation}

Let $P_{1,j}$ denote the spectral projection of $M_j$ for the eigenvalue $1$ and define  
\be\label{defpsi}
\psi_j:= P^*_{1,j}\psi_\cS, \ \ \ \ P_j:=|\psi_\cS\ket\bra \psi_j|.
\ee
Note that $\bra \psi_j|\psi_\cS\ket=1$ so that $P_j$ is a projection and, moreover, $M_j^*\psi_j=\psi_j$.  We introduce the following decomposition of $M_j$
\be\label{struct}
M_j:=P_j+Q_jM_jQ_j,  \ \ \ \mbox{with} \ \ \ Q_j=\one-P_j.
\ee
We denote the part of $M_j$ in $Q_j\C^d$, by 
$M_{Q_j}:=Q_jM_jQ_j$. 
It easily follows from these definitions that
\bea
P_jP_k=P_k  & & Q_jQ_k=Q_j \label{eq:ppqq}\\
Q_jP_k=0 & & P_kQ_j=P_k-P_j=Q_j-Q_k. \label{eq:pqqp}
\eea

\begin{prop}
\label{prop:psiform}
For any $n$, 
\be\label{eq:psi}
\Psi_n=|\psi_\cS\ket\bra\theta_n|+M_{Q_1}\cdots \mqn,
\ee
where
\bea
\theta_n & = & \psi_n+M_{Q_n}^*\psi_{n-1}+\cdots +M_{Q_n}^*\cdots M_{Q_2}^* \psi_1 \label{eq:theta1}\\
 & = & M_n^*\cdots M_2^* \psi_1 
\label{eq:theta2}
\eea
and where $\bra\psi_\cS,\theta_n\ket=1$.
\end{prop}

\proof Inserting the decomposition (\ref{struct}) into (\ref{def:psi}), and using \fer{eq:ppqq}, \fer{eq:pqqp}, we have
\be
\Psi_n= \sum_{j=1}^n P_j M_{Q_{j+1}}\cdots M_{Q_n} +M_{Q_1}\cdots \mqn.\nonumber
\ee
Since $P_j=|\psi_\cS\ket\bra \psi_j|,$ this proves (\ref{eq:psi}) and (\ref{eq:theta1}). From (\ref{eq:pqqp}), we obtain for any $j,k$,
\be
\label{mqm}
M_{Q_j}M_{Q_k}=M_{Q_j}M_k=Q_jM_jM_k.
\ee
Hence, $\Psi_n = P_1 M_1\cdots M_n +Q_1 M_1\cdots M_n = |\psi_\cS\ket\bra M_n^*\cdots M_2^* \psi_1| + M_{Q_1}\cdots \mqn$, 
which proves (\ref{eq:theta2}).
\qed

The following lemma is useful in the study of the asymptotic behaviour of $\Psi_n$. 
It is proven as the passage from \fer{eq:theta1} to \fer{eq:theta2} is, simply by considering the sequence $(M_{j_n})_n$ instead of $(M_n)_n$.

\begin{lem}
\label{lem:mmqrel} 
We have $\sum_{k=1}^n M_{Q_{j_1}}^*\cdots M_{Q_{j_k}}^*\psi_{j_{k+1}}=M_{j_1}^*\cdots M_{j_{n}}^*\psi_{j_{n+1}}$, for any sequence of indices $(j_n)_n$. 
\end{lem}


\subsection{Uniform bounds}\label{ssec:unifbound}

The operators $M_j$, and hence the product $\Psi_n$, are contractions on $\C^d$ for the norm $|||\cdot |||$.  In order to study their asymptotic
behaviour, we need some uniform bounds on the $P_j, Q_j,\ldots$ We denote by $\|\cdot\|$ the Euclidean norm on $\C^d$. The operator norm of a rank one projection $P=|\phi\ket\bra\chi|$ is simply $\|P\|=\|\phi\|\,\|\chi\|$ (see also after \fer{m13}). Recall also that $\|\psi_\cS\|=1$.

\begin{lem}\label{ubp}
We have $\ds \sup\, \{ \|M_{j_n}M_{j_{n-1}}\cdots M_{j_1}\|, \ n\in\N^*, \, j_k\in{\mathbb N}^*\} =C_0<\infty$. In particular, $\|\Psi_n\|\leq C_0$.
\end{lem}

\proof This follows from $|||M_j|||\leq 1$, $j\in\N^*$, and from the equivalence of the norms $\|\cdot\|$ and $|||\cdot|||$.
\qed

As a consequence, we get the following fundamental uniform estimates.
\begin{prop}\label{prop:unifbound} Let $C_0$ be as in Lemma \ref{ubp}. Then, the following bounds hold
\begin{enumerate}
\item For any $j\in\N^*$, $\|P_j\|=\|\psi_j\| \leq C_0$ and $\|Q_j\|\leq 1+C_0$.
\item $\ds \sup\, \{\|M_{Q_{j_n}}M_{Q_{j_{n-1}}}\cdots M_{Q_{j_1}}\|, \ n\in\N^*,\, j_k\in{\mathbb N}^*\} \leq C_0(1+C_0)$.
\item For any $n\in\N^*$, $\|\theta_n\|\leq C_0^2$.
\end{enumerate}
\end{prop}

\proof $1.$ Let $j\in\N^*$. By means of Von Neumann's ergodic Theorem, we have
$$
P_{1,j}=\lim_{N\ra\infty}\frac{1}{N}\sum_{k=0}^{N-1}M_j^k.
$$
Hence, we get from  Lemma \ref{ubp} that $\|P_{1,j}\|\leq C_0$. Since $\|\psi_\cS\|=1$ we also have $\|\psi_j\|\leq C_0$, and hence $\|P_j\|=\|\psi_j\| \leq C_0$ and $\|Q_j\|\leq 1+C_0$.

$2.$ Using (\ref{mqm}), we have 
$M_{Q_{j_n}}M_{Q_{j_{n-1}}}\cdots M_{Q_{j_1}}=Q_{j_n}M_{j_n}M_{j_{n-1}}\cdots M_{j_1},$
so that
$$
\|M_{Q_{j_n}}M_{Q_{j_{n-1}}}\cdots M_{Q_{j_1}}\|\leq \|Q_{j_n}\|C_0\leq C_0(1+C_0).
$$

$3.$ From (\ref{eq:theta2}) and the above estimates, we obtain $\|\theta_n\|\leq C_0\|\psi_1\|\leq C_0^2$.
\qed


\subsection{Asymptotic behaviour}\label{ssec:asympt}

We now turn to the study of the asymptotic behaviour of $\Psi_n$. If $M_n\equiv M\in\cme$ (recall the definition of $\cme$ given before Theorem \ref{thm1}), then $M_Q^n$ converges to zero, and $\Psi_n$ converges to the rank one projection $P=|\psi_\cS\ket\bra \psi|$ (exponentially fast). The following results shows that in the general case, $\Psi_n$ converges if and only if $\psi_n$ does.

\begin{prop}
\label{limpsi} 
Suppose that $\lim_{n\to\infty} \sup \{\|M_{Q_{j_n}}\cdots M_{Q_{j_1}}\|, \ j_k\in\N^*\}=0$. Then $\theta_n$ converges if and only if $\psi_n$ does. If they exist, these two limits coincide, and thus   
$$
\lim_{n\ra\infty}\Psi_n=|\psi_\cS\ket\bra \psi_\infty|,
$$
where $\psi_\infty=\lim_{n\rightarrow\infty} \psi_n$. Moreover, $|\psi_\cS\ket\bra \psi_\infty|$ is a projection.
\end{prop}

In general thus, we cannot expect to obtain pointwise convergence of the $\theta_n$, but we have to consider the ergodic average of $\theta_n(\po)$ in Theorem \ref{thm2}. This is natural in view of the interpretation of our results in terms of dynamical systems (a fluctuating system does not converge). For instance, if, in the random setting, $M(\omega)=P(\omega)=|\psi_\cS\ket \bra \psi(\omega)|$,  
then $\Psi_n(\po)=P(\omega_n)$ converges iff $\psi(\omega_n)$ does (this is actually a particular case of Proposition \ref{limpsi}). For a general random vector 
$\psi(\omega)$, there is no reason to expect a limit, due to 
the fluctuations of the $\psi(\omega_j)$. Nevertheless, Birkhoff's 
Ergodic Theorem asserts that the ergodic limit exists almost surely \cite{la}, 
and that 
\be\label{serl}
\lim_{N\ra\infty }\frac1N\sum_{j=1}^{N}P(\omega_j)={\mathbb E}(P)=
|\psi_{\mathcal{S}}\ket \bra {\mathbb E}(\psi)|=|\psi_{\mathcal{S}}\ket \bra \,\overline\psi \,|, \ \ \ a.s., \nonumber
\ee
where the projection $|\psi_{\mathcal{S}}\ket \bra \,\overline\psi \,|$ is constant in $\omega$. Theorem \ref{thm2} gives an analogous result in the general case.

\medskip
\noindent
{\bf Proof of Proposition \ref{limpsi}.\ } Suppose $\lim \theta_n=\theta_\infty$ exists, and  
let $\phi_n:=\theta_n-\theta_\infty$. We have $M_n^*\theta_\infty =  M_n^*\theta_{n-1}-M_n^*\phi_{n-1}
  = \theta_n-M_n^*\phi_{n-1}
  =  \theta_\infty+\phi_n-M_n^*\phi_{n-1}$. 
Hence, for any $n$, $M_{Q_n}^*\theta_\infty=Q_n^*\theta_\infty+Q_n^*\phi_n-M_{Q_n}^*\phi_{n-1}=: Q_n^*\theta_\infty+\chi_n$, 
where, by Proposition \ref{prop:unifbound}, $\lim_{n\to\infty} \chi_n=0.$ For any $p\in\N$, we can write $(M_{Q_n}^*)^p\theta_\infty = (M_{Q_n}^*)^{p-1}\theta_\infty+(M_{Q_n}^*)^{p-1}\chi_n$, 
and thus we obtain, by induction, $(M_{Q_n}^*)^p\theta_\infty=Q_n^*\theta_\infty+\sum_{k=0}^{p-1} (M_{Q_n}^*)^k\chi_n$. By assumption, given any $\eps>0$, we can chose $p_0$ large enough so that, for any $n$, $\|(M_{Q_n}^*)^{p_0}\| \leq \epsilon$. Using Proposition \ref{prop:unifbound}, we thus get that for any $n$, 
\be
\| Q_n^*\theta_\infty\| \ \leq \epsilon\|\theta_\infty\| + p_0C_0(C_0+1)\|\chi_n\|.\nonumber
\ee
Since $\chi_n$ goes to zero and $\eps$ is arbitrary, this proves that $\lim_{n\to\infty} Q_n^*\theta_\infty=0$. 
Now, $\langle \psi_\cS,\theta_\infty\rangle=1$ because $\langle \psi_\cS,\theta_n\rangle=1$ for all $n$. Thus we have $\psi_n=P_n^*\theta_\infty=\theta_\infty-Q_n^*\theta_\infty.$ Therefore $\lim \psi_n$ exists and equals $\theta_\infty$.

Conversely, suppose $\lim \psi_n=\psi_\infty$ exists. Since $M_n^*\psi_n=\psi_n$ for all $n$, we have 
$$
M_n^*\psi_\infty=M_n^*(\psi_\infty-\psi_n)+\psi_n.
$$
Combining this with Lemma \ref{ubp} shows that $\lim M_n^*\psi_\infty=\psi_\infty$. Furthermore, we have for any $n,p\in \N$, 
\bea
\theta_{n+p} & = & \sum_{j=1}^{n+p} M_{Q_{n+p}}^*\cdots M_{Q_{j+1}}^* \psi_j \nonumber \\
 & = & M_{Q_{n+p}}^*\cdots M_{Q_{n+1}}^*\sum_{j=1}^{n} M_{Q_n}^*\cdots M_{Q_{j+1}}^* \psi_j +  \sum_{j=n+1}^{n+p} M_{Q_{n+p}}^*\cdots M_{Q_{j+1}}^* \psi_j \nonumber \\
 & = & M_{Q_{n+p}}^*\cdots M_{Q_{n+1}}^* \theta_n + M_{n+p}^*\cdots M_{n+1}^* \psi_{n+1},\nonumber
\eea
where we used Lemma \ref{lem:mmqrel} to obtain the last line. It follows that 
\begin{eqnarray*}
\|\theta_{n+p}-\psi_\infty\| & \leq & \| M_{Q_{n+p}}^*\cdots M_{Q_{n+1}}^*\| \ \|\theta_n \|+ \|M_{n+p}^*\cdots M_{n+1}^* (\psi_{n+1}-\psi_\infty)\| \\
 & & + \|M_{n+p}^*\cdots M_{n+1}^* \psi_\infty-\psi_\infty\|
\end{eqnarray*}
Let $\eps>0$. By assumption, there exists $p_0$ such that for any $n$, $\| M_{Q_{n+p_0}}^*\cdots M_{Q_{n+1}}^*\| <\eps$. Thus we have, for all $n$,
\begin{eqnarray*}
\|\theta_{n+p_0}-\psi_\infty\| & \leq & \eps \|\theta_n\| +C_0 \|\psi_{n+1}-\psi_\infty \| \\
 & & + \Big\| M_{n+p_0}^*\psi_\infty-\psi_\infty+ \sum_{j=n+1}^{n+p_0-1}  M_{n+p_0}^*\cdots M_{j+1}^* (M_j^*\psi_\infty -\psi_\infty) \Big\| \\ 
 & \leq & \eps C_0^2 +C_0 \|\psi_{n+1}-\psi_\infty \| + \|M_{n+p_0}^*\psi_\infty-\psi_\infty\| \\
 & & +C_0(p_0-1) \sup_{n+1\leq j\leq n+p_0-1} \| M_j^*\psi_\infty-\psi_\infty\|.
\end{eqnarray*}
Since $\lim M_n^*\psi_\infty=\psi_\infty$, and since $\eps$ is arbitrary, we obtain $\lim \theta_n=\psi_\infty$. 
To prove that $|\psi_\cS\ket\bra\psi_\infty|$ is a projection, it suffices to note that $\bra\psi_n,\psi_\cS\ket=1$ for all $n$.
\qed

The previous convergence result relies on the
(exponential) decay of the product of operators $M_Q$. 
One can show that such a decay holds, in a deterministic situation, provided there exists $\delta>0$ such that 
$\|M_{Q_j}\|\leq e^{-\delta}$,  $\forall j\in\mathbb N$. One may obtain this bound for specific models, see for instance Section \ref{sssec:time}. To satisfy this bound, it is not only required that the $M_j$
belong to $\cme$ for all $j$, but also that
\begin {enumerate}
\item we have some uniformity in the spectrum of the $M_{Q_j}$, i.e. there exists a $\gamma>0$, such that 
\be\label{ev}
\sigma(M_{j})\setminus\{1\}\subset \{z\ | \ |z|\leq e^{-\gamma}\}, \ \ \forall j\in\N, \nonumber
\ee
\item we can prove that the corresponding spectral projectors do not behave too badly so that we can
still bound $M_{Q_j}$ efficiently, i.e. $\|Q_j\| e^{-\gamma}\leq e^{-\delta}$, for all $j$.
\end {enumerate}
Such conditions look quite stringent. However, Theorem \ref{thm1} shows that in the random setting, a similar exponential decay holds under much weaker assumptions.


\section{Proofs of Theorems \ref{thm1} and \ref{thm2}}
\label{sec:random}

\subsection{The probabilistic setting}\label{ssec:probadesc}

Let $(\Omega,\cF,\p)$ be a probability space underlying the random process. 
Given a random variable $f(\omega),$ we denote by
$\overline{f}:={\mathbb E}(f)=\int f(\omega)\d\p(\omega)$ its expectation. We define the probability measure $\d\P$ on $\Omega_\ext:=\Omega^{\N^*}$ in a standard fashion by
\be
\d\P =\Pi_{j\geq 1}\d \p_j, \ \ \ \mbox{where } 
\ \ \ \d \p_j\equiv \d \p, \ \  \forall j\in \N^*.\nonumber
\ee 
We denote points in $\Omega_{\rm ext}$ by $\po$. For any $r \in \N^*$, and for any collection of measurable sets in 
$\Omega$, 
$\{A_j\}_{1\leq j\leq r}$, we have 
\be
\P((\omega_1,\omega_2, \cdots, \omega_r )\in (A_1, A_2, \cdots, A_r))=
\Pi_{1\leq j\leq r}\int_{A_j}{\d \p(\omega_j)}.\nonumber
\ee
We also define the shift $T:\Omega_\ext\rightarrow\Omega_\ext$ by
\be\label{def:shift}
(T \po)_j= \omega_{j+1}, \ \ \forall \ \po 
=(\omega_j)_{j\in \N}\in \Omega_\ext.
\ee
$T$ is an ergodic transformation of $\Omega_\ext$. Further, we assume that the distribution of the matrices $M$ is caracterized by a measurable map
\be\label{def:mprocess}
\Omega\ni \omega_1 \mapsto M(\omega_1)\in M_d({\mathbb C}).
\ee
We write sometimes $M(\omega)$ instead of 
$M(\omega_1)$. Hence, for any set $B\in  M_d({\mathbb C})$, 
\be
\p (M(\omega)\in B)=\p(M^{-1}(B))=\int_{M^{-1}(B)}\d\p(\omega),\nonumber
\ee
and similarly for other random variables.
With these notations we have for the i.i.d. process
\bea
\Psi_n(\po):&=&
M(\omega_1)M(\omega_2)\cdots M(\omega_n)\\ \nonumber
&=&M(T^0\omega)M(T^1\omega)\cdots M(T^{n-1}\omega),\nonumber
\eea
where $\po\in\Omega^{{\mathbb N}^*}$.

In the same way as in (\ref{defpsi}), we introduce the random variable $\psi(\omega_1)\in {\mathbb C}^d$ 
defined as
\be\label{defpsi2}
\psi(\omega_1):= P_1(\omega_1)^*\psi_\cS,
\ee
where $P_1(\omega_1)$ denotes the spectral projection of $M(\omega_1)$ for the eigenvalue $1$, and we decompose 
\be\label{structrand}
M(\omega_1):=|\psi_\cS \ket\bra \psi(\omega_1)|+M_Q(\omega_1)=P(\omega_1)+M_Q(\omega_1)
\ee 
as in (\ref{struct}). Note that $\psi(\omega_1)$ and $M_Q(\omega_1)$ define {\it bona fide} random variables: $\omega_1\mapsto P_1(\omega_1)$ is measurable since $\omega_1\mapsto M(\omega_1)$ is \cite{az}. Consider finally the process
\bea
\theta_n(\po)&=&
M^*(T^{n-1}\omega)M^*(T^{n-2}\omega)\cdots M^*(T\omega)\psi(\omega)\\ \nonumber
&=&\sum_{j=1}^n M_Q^*(\omega_n)M_Q^*(\omega_{n-1})\cdots M_Q^*(\omega_{j+1})
\psi(\omega_j).
\eea
Note that $\theta_n$ is a Markov process, since $\theta_{n+1}(\po)=M^*_{n+1}(\omega_{n+1})\theta_{n}(\po)$.


\subsection{Decay estimates}
\label{ssec:probahyp}

In a first step, we derive a bound similar to the one in \fer{m17}, but only for some ``good $\omega$''. Then we will prove \fer{m17} via the Borel-Cantelli Lemma. 
\begin{lem}
\label{lem:goodneighb} 
Suppose $\p \left( M(\omega)\in \cme \right)>0.$ Then there exists an $M_0\in\cme$ such that for any $\eps>0$, $\p\left(\|M(\omega)-M_0 \|<\eps \right)>0.$
\end{lem}
\proof Let $M$ be such that $|||M|||\leq 1$ and $M\psi_\cS=\psi_\cS$. We decompose $M$ as in (\ref{struct}), and then $M\in\cme$ if and
only if $\sigma (M_Q)\subset \{z\in \C \,\ \mbox{s.t.} \ \, |z| < 1
\}.$ Note that it is equivalent to sr$(M_Q)<1$, where sr denotes the
spectral radius. For $n\in\N$, set 
\begin{eqnarray*}
\lefteqn{\cM_{(n)}:=}\\
&& \Big\{ M\in M_d(\C) \ {\rm s.t.}\  |||M|||\leq 1, M\psi_\cS=\psi_\cS, \sigma (M_Q)\subset \{z\in \C \,\ \mbox{s.t.} \ \, |z| \leq 1-1/n \} \Big\}.
\end{eqnarray*} 
The sets $\cM_{(n)}\subset M_d(\C)$ are compact, and we have $\cme = \bigcup_n \cM_{(n)}$. 
Thus, since $\p \left( M(\omega)\in \cme \right)>0$, there exists $n_0$ such that $\p
\left(M(\omega)\in \cM_{(n_0)}\right)>0.$ The lemma then follows with
$M_0\in \cM_{(n_0)}$ by the following standard argument using the compactness of $\cM_{(n_0)}.$ 

We construct a sequence $M_k\in\cM_{(n_0)}$ as follows. For any
$k\in\N^*$ we cover $\cM_{(n_0)}$ with a finite number of balls with centres in $\cM_{(n_0)}$  and radius $1/k$. For each fixed $k$, since $\p \left(M(\omega)\in \cM_{(n_0)}\right)>0$, there is at least one of these balls which has a non-zero probability. Pick one of these balls and denote its centre by
$M_k$. Thus, $\forall k\in\N^*$, $\p \left(
  \|M(\omega)-M_k\|\leq \frac{1}{k} \right)>0$. Since $\cM_{(n_0)}$ is compact, the sequence $M_k$ converges to some $M_0\in \cM_{(n_0)}$ (up to taking a subsequence). Now, for any given 
$\eps>0$, take $k$ large enough so that  $1/k<\eps/2$ and $\|M_k-M_0\|<\eps/2$.
Then
$$
\p\left(\|M(\omega)-M_0\|<\eps\right) \geq \p\left(\|M(\omega)-M_k\|<\eps/2\right)>0. 
$$ 
\qed

\begin{lem}
\label{probak} 
Suppose $\p \left( M(\omega)\in \cme \right)>0.$ Then, 
there exists $\widetilde\Omega_1\subset \Omega$ such that $p(\Omega_1)>0$,
and there exist constants $\delta>0$ and $C>0$ such that
$\forall k\in \N^*$, 
$$
(\omega_1, \omega_2, \cdots \omega_k)\in \widetilde\Omega_1^{k}  \ \ \ \Rightarrow \ \ \ 
\| M_Q(\omega_1)M_Q(\omega_2)\cdots M_Q(\omega_k)\|\leq C\e^{-\delta k}.
$$
\end{lem}
{\bf Proof.} Let $M_0$ be as in Lemma \ref{lem:goodneighb} and decompose it as $M_0=|\psi_\cS\ket \bra \psi_0 |+M_{Q_0}$.  
Using the Riesz representation,
$$
P_1=\frac{-1}{2 \pi i }\int_{\Gamma} (M-z)^{-1} dz,
$$
for a small contour $\Gamma$ encircling the eigenvalue one, one gets that the map
$M\mapsto M_Q$ is continuous on $\{ M\in M_d(\C) \, {\rm s.t.} \, \|M-M_0\| < \eta \}$, for $\eta$ small enough. Thus Lemma \ref{lem:goodneighb} shows that $\p \left( \|M_Q(\omega)-M_{Q_0}\| < \eps \right)>0.$ 

Hence, for any $\eps>0$, there exists $\Omega_\eps\subset \Omega$ such that $\p(\Omega_\eps)>0$, and such that for any $\omega_j\in \Omega_\eps$, we can write
\be
M_Q(\omega_j)= M_{Q_0}+\Delta(\omega_j)=: M_{Q_0}+\Delta_j \nonumber
\ee
with $\|\Delta_j\| < \eps$. With these notations, we have
\bea
\label{eq:mqdecomp}
\lefteqn{M_Q(\omega_1)M_Q(\omega_2)\cdots M_Q(\omega_k)}\\
 & = & \sum_{j=0}^k\ \ \ \sum_{1\leq k_1<k_2<\cdots <k_j\leq k}\!\!\!\!\!\!\!  M_{Q_0}^{k_1-1}\Delta_{k_1}
 M_{Q_0}^{k_2-k_1-1}\Delta_{k_2}\cdots M_{Q_0}^{k_j-k_{j-1}-1} \Delta_{k_j} M_{Q_0}^{k-k_j}.\nonumber
\eea
Since $M_0\in \cme$, we know that there exist $\gamma>0$ and $C>0$ such that $\|M_{Q_0}^n\| \leq C\e^{-\gamma n}$, for any $n$. 
Inserting this estimate in (\ref{eq:mqdecomp}), we get for any $k$ and $(\omega_1,\cdots,\omega_k)\in \Omega_\eps^k$,
\bea
\|M_Q(\omega_1)M_Q(\omega_2)\cdots M_Q(\omega_k)\| 
& \leq & \sum_{j=0}^k\ \ \sum_{1\leq k_1<k_2<\cdots <k_j\leq k} C^{j+1}\e^{-\gamma(k-j)}\eps^j \nonumber \\
 & = & C\e^{-\gamma k}\sum_{j=0}^k
\begin{pmatrix}k\cr j\end{pmatrix}(C\e^{\gamma}\eps)^j
= C\e^{-\gamma k}(1+C\e^{\gamma}\eps)^k.\nonumber
\eea
We now choose $\eps$ small enough so that  
$\e^{-\gamma}(1+C\e^{\gamma}\eps)=\e^{-\delta}$ 
for some $\delta>0$, and the proof is complete.
\qed


\subsection{Proof of Theorem \ref{thm1}} 

Let $\widetilde\Omega_1$ be as in Lemma \ref{probak}, and consider a product of $L$ matrices 
$$
M_Q(\omega_{j+1})M_Q(\omega_{j+2})\cdots M_Q(\omega_{j+L}).
$$
The $L$ matrices belong to
the set $M_Q(\widetilde\Omega_1)$ with positive probability, and by Lemma \ref{probak}, we have $\| M_Q(\omega_{j+1})M_Q(\omega_{j+2})\cdots M_Q(\omega_{j+L})\| \leq  C\e^{-\delta L}$. 
In other words, for any $j\in\mathbb N$, 
\be\label{block}
\P (\| M_Q(\omega_{j+1})M_Q(\omega_{j+2})
\cdots M_Q(\omega_{j+L})\|\leq C\e^{-\delta L})=: \cP\geq \p(\widetilde\Omega_1)^L>0.
\ee
We consider now a sequence of $nL$ matrices and estimate the probability 
that 
\be\label{esm}
\| M_Q(\omega_{1})M_Q(\omega_{2})
\cdots M_Q(\omega_{Ln})\| \leq (C\e^{-\delta L})^m C^{m+1}.
\ee
We decompose the above product into $n$ blocks of length $L$. Estimate (\ref{esm}) holds in particular if among these $n$ blocks, there are at least $m$ blocks in which (\ref{block}) holds. Indeed, let $j_1,\ldots, j_m \in \{0,\cdots ,n-1\}$ be such that 
\be\label{eq:goodblock}
\| M_Q(\omega_{j_kL+1})M_Q(\omega_{j_kL+2})
\cdots M_Q(\omega_{j_kL+L})\|\leq C\e^{-\delta L}, \ \ \ k=1,\cdots ,m.
\ee
We then isolate the $m$ ``good blocks'' and write
\bea
\lefteqn{\| M_Q(\omega_{1})M_Q(\omega_{2}) \cdots M_Q(\omega_{Ln})\|} \nonumber \\
 & \leq & \prod_{k=1}^m \| M_Q(\omega_{j_kL+1})\cdots M_Q(\omega_{j_kL+L})\| \ \|M_Q(\omega_1)\cdots M_Q(\omega_{j_1L})\| \nonumber \\
 & & \times \prod_{k=1}^{m-1} \| M_Q(\omega_{j_kL+L+1}) \cdots M_Q(\omega_{j_{k+1}L}) \| \ \|M_Q(\omega_{j_mL+L+1})\cdots M_Q(\omega_{nL})\|. \nonumber
\eea
Using Proposition \ref{prop:unifbound} and (\ref{eq:goodblock}), we obtain (\ref{esm}) with $C\geq C_0(1+C_0)$. Hence,
\be\label{eq:probdecay}
\P \Big(\| M_Q(\omega_{1}) \cdots M_Q(\omega_{Ln})\| \leq (C\e^{-\delta L})^m C^{m+1}\Big )\geq \sum_{j= m}^n\begin{pmatrix} n\cr j \end{pmatrix}{\cal P}^j (1-{\cal P})^{n-j}.
\ee
We now choose $m=[\beta n]$, for some $\beta>0$ to be chosen later, and where $[\cdot]$ denotes the integer part.
For $L$ large enough we have the following bound, uniform in $\beta$, 
\begin{equation}
\label{eq:probacontrol}
(Ce^{-\delta L})^{[\beta n]} C^{[\beta n]+1} \leq  C \e^{\delta L}(C^{2\beta}\e^{-\delta\beta L})^n =\tilde{C}\e^{-\gamma \beta n},
\end{equation}
for some $\gamma>0$. For $0< \beta<{\cal P}$, the following classical tail estimate follows from Hoeffding's inequality
\cite{hoef}
\bea\label{eq:minproba}
&& \sum_{j= [\beta n]}^n\begin{pmatrix} n\cr j
\end{pmatrix}{\cal P}^j
(1-{\cal P})^{n-j}=1-\sum_{j= 0}^{[\beta n]-1}\begin{pmatrix} n\cr j
\end{pmatrix}{\cal P}^j
(1-{\cal P})^{n-j}\nonumber\\
&& \qquad \qquad \qquad \qquad  \qquad \qquad \geq 1-e^{-2(n{\cal P}-[\beta n])^2/n}\geq 1-\e^{-\sigma n},
\eea
for $\sigma=2({\cal P}-\beta)^2$. Let us denote 
\be
G_n:=\Big\{\po=\{\omega_j\}_{j\geq 1} \ \mbox{\ s.t.\ } \ \| M_Q(\omega_{1})M_Q(\omega_{2})
\cdots M_Q(\omega_{Ln})\|\leq \tilde{C}\e^{-\gamma \beta n}\Big\}.\nonumber
\ee
Then, (\ref{eq:probdecay})-(\ref{eq:minproba}) show that, for all $n$, we have ${\mathbb P}(G_n)\geq 1-\e^{-\sigma n}$. 
Let $G_n^c$ be the complement of $G_n$ in $\Omega_\ext$. The first Borel-Cantelli Lemma
asserts that 
\be\label{bcl}
\sum_{n\geq 0} \P(G_n^c)<\infty \Rightarrow \P\big(\cap_n\cup_{m\geq n}G_m^c\big)=0.\nonumber
\ee
Since 
\be\label{exprob}
\P(G_n^c)\leq \e^{-\sigma n},
\ee 
we get that $G_n^c$ occurs
infinitely often with probability zero. Therefore, there exists
$\Omega_1\subset \Omega_\ext$, of measure one, such that for any $\po\in\Omega_1$,
there exists $n(\po)$ s.t. 
\be
\| M_Q(\omega_{1})M_Q(\omega_{2}) \cdots M_Q(\omega_{Ln})\| \leq
\tilde{C}\e^{-\alpha Ln}, \qquad \mbox{\ if\ $n\geq n(\po)$}
\label{nomega}
\ee 
where $\alpha=\gamma\beta/L$. The result of Theorem \ref{thm1} now follows: fix $L$ large enough and write $n=mL +p$ (with $0\leq p<L$). We conclude from Proposition \ref{prop:unifbound} and \fer{nomega} that if $\po\in \Omega_1$ (of measure one) and if $n\geq n_0(\po):=Ln(\po)+L$ (so that $m\geq n(\po)$), then
$$
\| M_Q(\omega_1)\cdots M_Q(\omega_n)\| \leq C_0\tilde C \e^{-\alpha L m}\leq C_0\tilde C \e^{\alpha p}\e^{-\alpha n}\leq C_0\tilde C \e^{\alpha L}\e^{-\alpha n}.
$$
This completes the proof of Theorem \ref{thm1}. \qed


\noindent {\it Remarks.} 1. We can take
$
n_0(\po):=\min\{m\geq 0 \, |\, \forall k\geq m, \| M_Q(\omega_{1})M_Q(\omega_{2})
\cdots M_Q(\omega_{k})\|\leq \tilde{C}\e^{-\alpha k}\},
$
so that, defining $g_n=\Big\{\po\  \mbox{s.t.}\,   \| M_Q(\omega_{1})M_Q(\omega_{2})
\cdots M_Q(\omega_{n})\|\leq \tilde{C}\e^{-\alpha n}\Big\}$, we have
$
\Big\{\po \ |  \ n_0(\po)=m \Big\}=\cap_{k\geq m}g_k\cap g^c_{m-1}.
$
Hence, by the estimate (\ref{exprob})  and our definition of $\alpha, \beta$
and $\sigma$, we get  
\be\label{ddd}
{\mathbb E}(e^{\alpha n_0})\leq \sum_{m\geq 0}e^{\alpha m}
{\mathbb P}(g^c_{m-1}) \leq \e^{\sigma/L} \sum_{m\geq 0}e^{(\gamma\beta-\sigma)m/L}= \e^{\sigma/L}
 \sum_{m\geq 0}e^{(\gamma\beta-2({\cal P}-\beta)^2)m/L} <\infty,
\ee
for $\beta$ small enough.

2. By Proposition \ref{prop:unifbound}, \fer{m17} holds for all
$\po\in\Omega_1$ and
for all $n\in\N^*$, up to replacing $\tilde{C}$ by
$C_{\po}=\tilde{C}e^{\alpha n_0(\po)}$.\\

As we mentioned in the introduction, if one can chose
$\psi(\omega)\equiv\psi$ one may expect a convergence without any
ergodic averaging. Indeed, the following result is a direct
consequence of (\ref{eq:theta2}) and Theorem \ref{thm1}  


\begin{prop} Let $M(\omega)$ be a random reduced dynamics
  operator. Suppose that $\p(M(\omega)\in\cme)>0$ and there exists
  $\psi\in \C^d$ such that $M(\omega)^*\psi=\psi$ for all $\omega$. Then there exists a set $\Omega_1\subset\Omega^{\N^*}$, and a constant $\alpha>0$, s.t. ${\mathbb P}(\Omega_1)=1$ and s.t. for any $\po\in\Omega_1$, there exists $C_{\po}$ so that for any $n\in{\N^*}$,
\begin{equation}
\left\| M_1(\omega)\cdots M_n(\omega) - \frac{|\psi_\cS\ket \bra \psi|}{\bra \psi,\psi_\cS\ket} \right\|\leq C_{\po}\e^{-\alpha n}.\nonumber
\end{equation}
\end{prop}

This result applies in particular to bistochastic matrices with
$\psi_\cS=\psi=\frac{1}{\sqrt{d}}(1,\cdots,1).$ Infinite product of doubly stochastic matrices have been considered in e.g. \cite{sch}, where the author looks for sufficient conditions which ensure that such a product converges to the rank one projection $|\psi_\cS\ket\bra \psi_\cS|.$ However, only deterministic products are considered.


\subsection{Proof of Theorem \ref{thm2}}\label{ssec:ergocvg}

The proof of Theorem \ref{thm2} relies on the estimate
\fer{m17}, but not directly on the more stringent assumption
$\p\left(M(\omega)\in\cme\right)>0$. Consider the condition  


\begin{itemize}
\item[{\bf (RE)}] There exist $\tilde\Omega\subset\Omega^{{\mathbb N}^*}$, $C,\alpha>0$ such that ${\mathbb P} (\tilde\Omega)=1$, and such that for any 
$\po\in \tilde\Omega$, there exists  $n_0(\po)\in\mathbb N$ so that  $\|M_{Q}(\omega_1)\cdots M_{Q}(\omega_n)\|\leq Ce^{-\alpha n}$ if $n\geq n_0(\po)$, and, $\|M_{Q}(\omega_1)\cdots M_{Q}(\omega_n)\|\leq C$, otherwise.
\end{itemize}
Theorem \ref{thm1} and Proposition \ref{prop:unifbound} shows that (RE) is a consequence of 
 $\p\left(M(\omega)\in\cme\right)>0$. Consequently, Theorem \ref{thm2} follows from Theorem \ref{thm:ergodic} below. We denote by $\ds P_{1,\overline{M}}$ the spectral projection of
$\overline{M}:=\E(M(\omega))$ for the eigenvalue $1$. 

\noindent{\it Remark.\ } Note that
$M(\omega)\psi_\cS=\psi_\cS$ for all $\omega$ implies $\overline{M}\psi_\cS=\psi_\cS.$
In general, it may happen that $1$ is a degenerate eigenvalue for $\overline{M}$ even though it is non-degenerate for any $M(\omega)$. This will however not be the case here as the following lemma shows.  
\begin{lem}\label{lem:mbar} Suppose Assumption {\bf (RE)} holds, then
  $\overline{M}\in\cme$. Moreover, we have
  $P_{1,\overline{M}}=|\psi_\cS \ket\bra \theta|$,
where $\theta$ satisfies 
\be\label{ttb}
\theta = (\one-\overline{M_Q^*})^{-1}\overline{\psi} = \sum_{k=0}^\infty \overline{M_Q^*}^{k}\overline{\psi} = P_{1,\overline{M}}^*\overline\psi = P_{1,\overline{M}}^*\psi_\cS. 
\ee
\end{lem}
\medskip
\noindent \proof We have $\overline{M}=|\psi_\cS\ket\bra \overline \psi|
+\overline{M_Q}.$ Since $M(\omega)\psi_\cS=\psi_\cS$ and $\bra
\psi(\omega)|\psi_\cS\ket=1$ for any $\omega$, we get
$\overline{M_Q}\psi_\cS=0$ and $\bra \overline\psi
|\psi_\cS\ket=1$. Thus, for any $n\in\N$ we have (see also (\ref{eq:psi})-(\ref{eq:theta1}))
\be\label{eq:mbarstruct}
\overline{M}^n=|\, \psi_\cS \, \ket \bra \, \overline\psi+\overline{M_Q}^*\overline\psi+\cdots +\overline{M_Q^*}^{n-1}\overline\psi \,| + \overline{M_Q}^n.
\ee
By independence, we
have, for any $k\in\mathbb N,$
$\overline{M_Q}^{k}=\int \d{\mathbb P}M_{Q}(\omega_1)\cdots  M_{Q}(\omega_k)$
so that
\be
\|\overline{M_Q}^{k}\|\leq \int \d{\mathbb P}\|M_{Q}(\omega_1)\cdots  M_{Q}(\omega_k)\|.
\ee
Due to Assumption {\bf (RE)} and considerations similar to (\ref{ddd}), we can bound the right hand side by
\be
{\mathbb P}\{n_0(\po)>k\}C+{\mathbb P}\{n_0(\po)\leq k\}Ce^{-\alpha k}\leq
\bar C e^{-\alpha k},
\ee
for some constant $\bar C$.
Hence we have
\be\label{mm130}
\|\overline{M_Q}^{k}\|\leq Ce^{-\alpha k}, \ \ \forall k\in\N.
\ee
In particular, the spectrum of $\overline{M_Q}$ lies inside the open complex unit disk.
 
This estimate and (\ref{eq:mbarstruct}) prove $\ds \lim_{n\to\infty} \overline{M}^n= |\psi_\cS \ket\bra \theta|$ where 
$$
\theta=\sum_{k=0}^\infty \overline{M_Q^*}^{k}\overline{\psi}=(\one-\overline{M_Q^*})^{-1}\overline{\psi}.
$$
Hence $\overline{M}\in\cme$ and $P_{1,\overline{M}}=|\psi_\cS \ket\bra\theta|$.
The last identities in (\ref{ttb}) are immediate since $\bra \psi_\cS,\overline\psi \ket=\bra\psi_\cS,\psi_\cS\ket=1$.
\qed

\noindent {\it Remarks.\ } 1. In general $\overline\psi$ is not an eigenvector of $\overline{M^*}$, and $\theta\neq\overline\psi$; see Theorem \ref{thm:spinatom} and the remark thereafter.

2. Without Assumption {\bf (RE)}, instead of \fer{mm130}, we can only prove that $\|\overline{M_Q}^{k}\|\leq C_0(1+C_0)$  and that $\|\overline{M}^{k}\|\leq C_0$, $\forall k\in\mathbb N$, which is not enough to exclude that $1\in \sigma (\overline{M_Q})$.

\begin{thm}
\label{thm:ergodic}
Suppose Assumption {\bf (RE)} holds. Then, there exists a set $\Omega_0\subset\Omega^{\N^*}$ such that $\P(\Omega_0)=1$ and, for all $\po\in \Omega_0$, 
\be
\lim_{N\ra\infty}\frac1N\sum_{n=1}^N \theta_n(\po)= \theta.\nonumber
\ee
\end{thm}

\noindent
{\bf Proof of Theorem \ref{thm:ergodic}. } Recall that $T$ denotes the shift map on $\Omega_\ext$ (see (\ref{def:shift})). 
Using (\ref{eq:theta1}), we get
\bea\label{eq:thetaerg}
\sum_{n=1}^N\theta_{n}(\po) & = & \sum_{n=1}^N \sum_{j=0}^{n-1} M_Q^*(T^{n-1}\omega)\cdots M_Q^*(T^{j+1}\omega)\psi(T^j\omega) \nonumber\\
 & = & \sum_{k=1}^N \sum_{j=0}^{N-k} M_Q^*(T^{k+j-1}\omega)\cdots M_Q^*(T^{j+1}\omega)\psi(T^j\omega).
\eea
Let us introduce random vectors $\theta^{(k)}$, $k=1, 2, \ldots,$  by
\bea\label{thek}
\theta^{(k)}(\po) & = & \theta^{(k)}(\omega_1,\cdots,\omega_k)\\
&=& M_Q^*(T^{k-1}\omega)M_Q^*(T^{k-2}\omega)\cdots M_Q^*(T\omega)\psi(T^0\omega)
\nonumber\\
&=&M_Q^*(\omega_k)\theta^{(k-1)}(\omega_1,\cdots,\omega_{k-1}),
\nonumber
\eea
so that, by (\ref{eq:thetaerg}),
\bea\label{eq:thetaerg2}
\frac{1}{N}\sum_{n=1}^N\theta_{n}(\po)&=&\sum_{k=1}^N\frac{1}{N}\sum_{j=0}^{N-k}
\theta^{(k)}(T^j\po)\nonumber\\
&=&\sum_{k=1}^\infty
\chi_{\{k\leq N\}}\sum_{j=0}^{N-k}\theta^{(k)}(T^j\po)\frac{1}{N}
=: \sum_{k=1}^\infty g(k,N, \po).
\nonumber
\eea
For each fixed $k$, by ergodicity, there exists a set $\Omega_{(k)}\subset\Omega^{\N^*}$ of probability one, such that, for all $\omega\in \Omega_{(k)}$, the following
limit exists
\bea
\lim_{N\ra\infty}g(k,N, \po)&=&\lim_{N\ra\infty}
\frac{1}{N-k+1}\sum_{j=0}^{N-k}\theta^{(k)}(T^j\po)\frac{N-k+1}{N}\nonumber\\
&=&\lim_{M\ra\infty}\frac{1}{M+1}\sum_{j=0}^{M}\theta^{(k)}(T^j\po)
={\mathbb E}(\theta^{(k)}).\nonumber
\eea
Therefore, on the set $\Omega_{\infty}:=\cap_{k\in\mathbb N}\Omega_{(k)}$ of 
probability one,  for any $k\in\mathbb N$, we have by independence of the 
$M(\omega_j)$,
\be\label{eq:limg}
\lim_{N\ra\infty} g(k,N, \po)={\mathbb E}(\theta^{(k)})=
(\overline{M_Q^*})^{k-1}\overline{\psi}.
\ee
Also, by Assumption {\bf (RE)} and Proposition \ref{prop:unifbound}, on $\Omega_0=\tilde\Omega \cap\Omega_{\infty}$, 
a set of probability $1$, 
we have the bound $\|\theta^{(k)}(T^j\po)\|\leq C_0Ce^{\alpha n_0(T^j\po)} e^{-\alpha(k-1)}$. Hence, for all $N$ large enough, and for all $1\leq k\leq N$, by ergodicity
\bea
\|g(k,N, \po)\|&\leq&  C_0C \frac{1}{N}\sum_{j=0}^{N-1}e^{\alpha n_0(T^j\po)}e^{-\alpha(k-1)}\leq
2C_0C {\mathbb E}(e^{\alpha n_0}) e^{-\alpha(k-1)},
 \nonumber
\eea
which is finite by (\ref{ddd}). For $k>N$, $\|g(k,N, \po)\|=0\leq 2C_0C{\mathbb E}(e^{\alpha n_0}) e^{-\alpha(k-1)}$. 
Since $e^{-\alpha(k-1)}\in l^1(\mathbb N)$,
we can apply the Lebesgue Dominated Convergence Theorem in (\ref{eq:thetaerg2}) to conclude
that $\lim_{N\to\infty} \frac{1}{N}\sum_{n=1}^N
\theta_n(\po)=\theta$ a.s. (recall  Lemma \ref{lem:mbar}). 
\qed


\subsection{The Mean Ergodic Theorem and Lyapunov exponents}\label{ssec:met}

We present in this section some results for products ``in reverse order'' of the form 
$\Phi_n(\po):= M(\omega_n) \cdots M(\omega_1)$, as well as results related to the Multiplicative Ergodic Theorem. The following results are standard, see  e.g. \cite{la}, Theorem 3.4.1 and Remark 3.4.10 ii). The limits 
$$
\Lambda_\Phi(\po)=\lim_{n\ra\infty}(\Phi_n(\po)^*\Phi_n(\po))^{1/2n} \ \  {\rm\  and\ } \ \ 
\Lambda_\Psi(\po)=\lim_{n\ra\infty}(\Psi_n(\po)^*\Psi_n(\po))^{1/2n}
$$
exist almost surely, the top Lyapunov exponent of $\Lambda_\Phi(\po)$ and of $\Lambda_\Psi(\po)$ are the same, called $\gamma_1(\po)$; $\gamma_1(\po)$ is constant a.s., and so is its multiplicity. It is usually difficult to prove
that the multiplicity of $\gamma_1(\po)$ is $1$. In the present case however, we can get this result.
\begin{thm}
\label{prop:randlyap} 
Suppose $\p(M(\omega)\in\cme)>0$. Then there exist $\alpha>0$, a random vector $\eta_\infty(\po)$ and $\Omega_0\subset\Omega^{{\mathbb N}^*}$ with $\P(\Omega_0)=1$ such that for any $\po\in\Omega_0$ and $n\in\N^*$
\be\label{phicvg}
\Big\| \Phi_n(\po)- |\psi_\cS\ket\bra\eta_\infty(\po)| \Big\| \leq C_{\po} e^{-\alpha n}, \ \ \mbox{for some } \  C_{\po}.
\ee
As a consequence, for any $\po\in\Omega_0$, $\gamma_1(\po)$ is of multiplicity one. 
\end{thm}

\proof We first decompose $\Phi_n(\po)$ in the same way as $\Psi_n(\po)$. We have, 
\be\label{eq:phi}
\Phi_n(\po)=|\psi_\cS\ket\bra\eta_n(\po)|+M_Q(\omega_n)\cdots M_Q(\omega_1),
\ee
where
\bea
\eta_n(\po) & = & \psi(\omega_1)+M_Q(\omega_1)^*\psi(\omega_2)+\cdots +M_Q(\omega_1)^*\cdots M_Q(\omega_{n-1})^* \psi(\omega_n) \ \ \ \label{eq:eta1}\\
 & = & M(\omega_1)^*\cdots M(\omega_n)^* \psi(\omega_n) \label{eq:eta2}
\eea
and $\bra\psi_\cS,\eta_n(\po)\ket=1$. Using (\ref{eq:eta1}) and Lemma \ref{lem:mmqrel}, we have for any $n,p\in \N$
\begin{eqnarray*}
\lefteqn{
\eta_{n+1+p}(\po)-\eta_{n+1}(\po)}\\
 & = & M_Q(\omega_1)^*\cdots M_Q(\omega_n)^*\sum_{j=1}^p
 M_Q(\omega_{n+1})^*\cdots M_Q(\omega_{n+j})^*\psi(\omega_{n+j+1})\nonumber\\
 & = & M_Q(\omega_1)^*\cdots M_Q(\omega_n)^* M(\omega_{n+1})^*\cdots M(\omega_{n+p})^* \psi(\omega_{n+p+1}).
\end{eqnarray*}
Thus, Lemma \ref{ubp} and Proposition \ref{prop:unifbound} imply
\be
\|\eta_{n+p+1}(\po)-\eta_{n+1}(\po)\| \leq C_0^2\|M_Q(\omega_{n})M_Q(\omega_{n-1})\cdots M_Q(\omega_1)\|,\nonumber
\ee
uniformly in $p\geq 0$. Using Theorem \ref{thm1}, there exists
$\Omega_0$ with $\P(\Omega_0)=1$ and $\alpha>0$ such that for any $n,p\in\N$, $\po\in \Omega_0$, there exists $C_{\po}$ s.t. $\|\eta_{n+p+1}(\po)-\eta_{n+1}(\po)\|\leq C_0^2C_{\po} \e^{-n\alpha}$. 
Therefore, $\eta_n(\po)=\eta_\infty(\po)+O_{\po}(\e^{-n\alpha})$. Together with (\ref{eq:phi}) and Theorem \ref{thm1}, this proves (\ref{phicvg}).

Now, we can write
\bea
\Phi_n(\po)^*\Phi_n(\po)&=&|\eta_\infty(\po) \ket\bra \eta_\infty(\po) |+ O_{\po}(e^{-n\alpha})\nonumber\\
  &=&\|\eta_\infty(\po)\|^2|\hat\eta(\po) \ket\bra\hat \eta(\po) |+ O_{\po}(e^{-n\alpha}),\nonumber
\eea
where $\hat\eta(\po)=\eta_\infty(\po)/\|\eta_\infty(\po)\|$, so that $|\hat\eta(\po)\ket\bra
\hat\eta(\po)|$ is an orthogonal projector. Hence, by perturbation theory, the non-negative self-adjoint operator
$\Phi_n(\po)^*\Phi_n(\po)$ has a simple eigenvalue 
$\|\eta_\infty(\po)\|^2+O_{\po}(e^{-n\alpha})$ and all its other 
eigenvalues are non-negative and bounded from above by $C_{\po}e^{-n\alpha}$, for some $C_{\po}>0$. Thus the
$2n$'th root of $\Phi^*_n(\po)\Phi_n(\po)$ has in the limit $n\ra\infty$, 
a simple eigenvalue 1, and 
all its other eigenvalues are bounded from above by $e^{-\alpha/2}$.
\qed


\section{From repeated interactions to matrix products}
\label{sec:reduction}

In this section, we show how to reduce the evolution of states of a
R.I. quantum system (for observables on the small system 
${\mathcal S}$) to the study of a product of non-identical matrices
given by a RDO process.


\subsection{The repeated interaction model}\label{ssec:rimodel}

We call the {\it time step} $m$ the index $m\geq 1$ that labels the interaction step during which the system ${\mathcal S}$
interacts with the element ${\mathcal E}_m$ of the chain 
${\mathcal C}={\mathcal E}_1+{\mathcal E}_2+{\mathcal E}_3+\cdots$, eventhough the duration $\tau_m$ of this
interaction may vary with $m$, with the understanding that $\inf_{m\in\N^*}\tau_m>0$.

The Hilbert space $\cH_\cS$ is finite-dimensional. We denote by $L_{\cE_k}$ the free Liouvillean of $\cE_k$, acting on $\cH_{\cE_k}$. The invariant reference vector 
$\psi_{\cE_k}\in {\cH_{\cE_k}}$ is cyclic and separating for
$\fm_{\cE_k}$, it is independent of time and satisfies $L_{\cE_k}\psi_{\cE_k}=0$.  The interaction operator between the $m$th element of the chain and
the reference system at time $m$ is a self-adjoint operator $V_m\in \fm_\cS\otimes \fm_{\cE_m}$. 
The coupled dynamics on $\cH_\cS\otimes{\mathcal H}_{\cE_m}$ is then generated by the interaction Liouvillean $L_m=L_\cS+L_{\cE_m}+V_m$. 
Moreover, if $J_m$ and $\Delta_m$ denote the modular data of the pair $( \fm_\cS\otimes \fm_{\cE_m}, \psi_\cS\otimes \psi_{\cE_m})$ (see e.g. \cite{BR}), we further require 
\be
\Delta_m^{1/2}V_m\Delta_m^{-1/2}\in \fm_\cS\otimes \fm_{\cE_m}, \ \ \forall m\geq 1.\nonumber
\ee
The evolution operator acting on ${\cal H}={\cal H}_\cS\otimes {\cal H}_{\cC}$ in the Schr\"odinger picture at time step $m$ is given by the unitary $U(m)=e^{-i \widetilde{L}_m}e^{-i \widetilde{L}_{m-1}}\cdots e^{-i \widetilde{L}_1}$, where $\widetilde{L}_m$ is given in \fer{m3}.

Let us finally define the Liouville operator $K_m$ at time $m$, see \cite{bjm,jp2002,MMS}, by
\be
K_m=\tau_m\left[ L_\cS+L_{\cE_m}+V_m-J_m\Delta_m^{1/2}V_m\Delta_m^{-1/2}J_m\right].\nonumber
\ee
Its main dynamical features, 
\begin{eqnarray}
e^{i\tau_mL_m} A e^{-i\tau_mL_m}&=&e^{iK_m} A e^{-iK_m}, \mbox{\ \ for $A\in{\frak M}_\cS\otimes{\frak M}_\cE$, $t\in{\mathbb R}$,}\label{cliouv}\\
K_m\psi_\cS\otimes\psi_{\cE,m} &=&0,
\label{cliouv2}
\end{eqnarray}
(see also \fer{m8},\fer{m8.1}), are proven to hold by using standard relations of the modular data $\Delta_m, J_m$, see e.g. \cite{bjm}. 
Note also the bound $\|e^{\pm iK_m}\|\leq e^{\tau_m \|\Delta_m^{1/2}V_m\Delta_m^{-1/2}\|}$.

We proceed along the same lines as in \cite{bjm}, taking into account that in the present setup, $V_m$ may vary with $m$. Let $\cH_\cC:= \otimes_{k\geq 1} \cH_{\cE_k}$,
where the infinite tensor product is understood with respect to the reference vector $\psi_\cC:=\otimes_{k\geq 1}\psi_{\cE_k}$. Since $L_m$ and $L_{\cE_k}$  
commute for $k\neq m$, we can write successively
\bea
e^{-i\widetilde{L}_1}&=&e^{-i\tau_1L_1}\ e^{-i\tau_1\sum_{k>1}L_{1_k}}, \nonumber\\
e^{-i\widetilde{L}_2}&=&e^{-i\tau_2 L_{\cE_1}}e^{-i\tau_2 L_2}\
e^{-i\tau_2\sum_{k>2}L_{2_k}},\nonumber\\
&\vdots&\nonumber\\
e^{-i\widetilde{L}_m}&=&e^{-i\tau_m\sum_{k<m}L_{\cE_k}}\ e^{-i\tau_mL_m} \
e^{-i\tau_m\sum_{k>m}L_{\cE_k}},\nonumber
\eea
so that we obtain
\bea
\label{sch}
\lefteqn{U(m)}\nonumber\\
&=&\exp\left[ -i \sum_{j=1}^m \sum_{k=1}^{j-1}\tau_jL_{\cE_k}\right] \ e^{-i \tau_m L_m}e^{-i\tau_{m-1} L_{m-1}}\cdots e^{-i\tau_1 L_1}
\ \exp\left[ -i \sum_{j=1}^m \sum_{k>j}\tau_jL_{\cE_k}\right]\nonumber\\
&\equiv&U_m^- \ e^{-i \tau_m L_m}e^{-i\tau_{m-1} L_{m-1}}\cdots e^{-i\tau_1 L_1}\ U_m^+.
\eea
The unitaries $U_m^\pm$ act trivially on $\cH_\cS$ and satisfy $U_m^\pm\psi_\cC=\psi_\cC$, $\forall m\in\N^*$. 
The corresponding Heisenberg evolution of an observable $O\in \fm=\fm_\cS\otimes\fm_\cC$ at time $m$ is \fer{m5.1}. 
Given a normal state $\varrho$ on $\fm$, we want to understand the limit of large $m$ of the expectation value $\varrho(\alpha_{RI}^m(O)) = \varrho(U(m)^*\ O\ U(m))$, 
for observables on the small system $\cS$, i.e. $O=A_\cS\otimes\one_\cC$ (more general observables, such as the so-called
instantaneous observables following the terminology of \cite{bjm}, will be considered in \cite{bjm2}).
In other words, we wish to adress the long time behaviour of the evolution of the state $\varrho$ defined by
$\varrho\circ \alpha_{RI}^m $ on $\fm$. 

Any normal state is a convex combination of vector states, i.e., its density matrix is given by a convex combination of rank-one projections,
$\rho = \sum_n p_n |\phi_n\rangle\langle\phi_n|$, 
where $\sum_np_n=1$, and where the $\phi_n$ are normalized vectors in $\cH$. It is then clear that convergence of the vector states $|\phi_n\rangle\langle\phi_n|$ implies convergence of $\rho$ (see \cite{bjm} for more details). Therefore, we will limit our attention to vector states $\scalprod{\phi}{\cdot\ \phi}$ from now on. Moreover, because $\psi:=\psi_\cS\otimes\psi_\cC$ is cyclic for the commutant $\fm'$, it is enough to consider $\phi$ of the form
\be
\label{psib}
\phi = B' \psi_\cS\otimes\psi_\cC=:B'\psi,
\ee
for some 
\begin{equation}\label{bpn}
B'= B_\cS'\otimes_{n=1}^N B_n'\otimes_{n>N}\one_{\cE_n}\ \in\fm',
\end{equation}
with $B_\cS'\in\fm_\cS'$, $B_n'\in\fm_{\cE_n}'$ (up to a vanishing 
error as $N\ra \infty$, see \cite{bjm}). Hence, we are lead to to consider 
\be\label{basic}
\bra \psi ,(B')^* \alpha_{RI}^m(O) B' \psi\ket=
\bra \psi , (B')^*  B'\alpha_{RI}^m(O) \psi\ket, 
\ee
for observables $O=A_\cS\otimes \one_\cC$. In \fer{basic} we use that $B'$ commutes with $\alpha_{RI}^m(O)$.


\subsection{Observables on the small system}\label{ssec:systobs}

In this paper we restrict our attention to the Heisenberg evolution of observables $A_\cS\in \fm_\cS$. That is, we consider
\bea
&&\alpha_{RI}^m(A_\cS\otimes \one)= U(m)^*\ (A_\cS\otimes \one)\ U(m)\nonumber\\
&& \quad \quad = (U_m^+)^* \ e^{i \tau_1{L_1}}e^{i\tau_2{L_{2}}}\cdots e^{i \tau_m{L_m}} \ (A_\cS\otimes \one )\
 e^{-i \tau_m{L_m}}e^{-i\tau_{m-1}{L_{m-1}}}\cdots e^{-i \tau_1{L_1}}\ U_m^+,\nonumber
\eea
where we made use of the fact that $U_m^-$ acts trivially on
$\cH_\cS$. 
Considering (\ref{basic}), we note that due to the properties of the unitary $U^{+}(m)$
and those of the Liouvillean, (\ref{cliouv})-(\ref{cliouv2}),  we have
\bea
\lefteqn{\alpha_{RI}^m(A_\cS\otimes \one)\, \psi}\nonumber\\
 &=&
U^{+}(m)^* \ e^{i \tau_1{L_1}}e^{i\tau_2{L_{2}}}\cdots e^{i \tau_m{L_m}} \ (A_\cS\otimes \one )\
 e^{-i \tau_m{L_m}}e^{-i\tau_{m-1}{L_{m-1}}}\cdots e^{-i \tau_1{L_1}} \psi\nonumber\\
&=&U^{+}(m)^* \ e^{i {K_1}}e^{i{K_{2}}}\cdots e^{i {K_m}} \ (A_\cS\otimes \one )\psi.\nonumber
\eea
Let us introduce $P_N=\one_\cS\otimes \one_{\cE_1}\otimes \cdots \one_{\cE_N}\otimes P_{\psi_{\cE_{N+1}}}\otimes 
P_{\psi_{\cE_{N+2}}} \otimes \cdots $, where 
$P_{\psi_{\cE_k}}=|\psi_{\cE_k}\ket\bra \psi_{\cE_k} |$. Due to (\ref{bpn}),
we have $\bra \psi | (B')^*  B'=\bra \psi | (B')^*  B'P_N$. Moreover, introducing the unitary operator
\be
\tilde U^{+}(N)=\exp\left[ -i \sum_{j=1}^{N-1} \sum_{k=j+1}^N \tau_j L_{\cE_k}\right]=U^{+}(m)P_N,\nonumber
\ee
which is {\it independent} of $m$, we can write for $m>N$,
\bea
\lefteqn{\bra \psi, (B')^*  B'\alpha_{RI}^m(A_\cS\otimes\bbbone) \psi\ket}\nonumber\\
 &=& \bra \psi, (B')^*  B' \tilde U^{+}(N)^*  P_N e^{i{K_1}}e^{i{K_{2}}}\cdots e^{i {K_m}}  (A_\cS\otimes \one ) \psi\ket\nonumber\\
 &=& \bra \psi, (B')^*  B' \tilde U^{+}(N)^*  e^{i {K_1}}\cdots e^{i{K_{N}}}P_Ne^{i{K_{N+1}}}\cdots e^{i {K_m}}  (A_\cS\otimes \one ) \psi\ket\nonumber.
\eea
At this level, introducing the projector $P=\one_\cS\otimes |\psi_\cC\ket\bra \psi_\cC|$, satisfying $P\psi=\psi$, we observe that $P_Ne^{i{K_{N+1}}}\cdots e^{i {K_m}}(A_\cS\otimes \one )\psi=Pe^{i{K_{N+1}}}\cdots e^{i {K_m}}P(A_\cS\otimes \un )\psi$.
Thus, we can invoke Proposition 4.1 of \cite{bjm} which states that for any $q\geq 1$ and any
distinct integers $n_1, \ldots, n_q$,
\be
P e^{i{K_{n_1}}}e^{i{K_{n_2}}}\cdots e^{i {K_{n_q}}}P=P e^{i {K_{n_1}}}Pe^{i{K_{n_2}}}P\cdots Pe^{i {K_{n_q}}}P.\label{mm100}
\ee
We note that the proof of this statement does not require that the $K_n$'s be identical as operators on $\cH_\cS\otimes \cH_\cE$. Therefore, introducing operators $M_j$ acting on $\cH_\cS$ by 
\be\label{defmj}
P e^{i{K_{j}}}P=M_j\otimes |\psi_\cC\ket\bra \psi_\cC|, \ \ \ \mbox{or} \ \ \ M_j\simeq P e^{i {K_{j}}}P,
\ee
we arrive at
\bea\label{almost}
&&\bra \psi | (B')^*  B'\alpha_{RI}^m(A_\cS\otimes\one) \psi\ket\\
&&\quad\quad=\bra \psi | (B')^*  B' \tilde U^{+}(N)^* e^{i {K_1}}\cdots 
e^{i{K_{N}}}P M_{N+1} M_{N+2} \cdots M_m  (A_\cS\otimes \one )
\psi\ket.\nonumber
\eea
Hence, given the repeated interactions dynamics, in order to understand the long time behaviour of the 
evolution of any states $\varrho\circ \alpha_{RI}^m$ on the set of observables on the reference system $\cS$, it is 
necessary (and sufficient) to address the asymptotic behaviour of the product of matrices $M_j\simeq P e^{i {K_{j}}}P$ of the type $\Psi_m=M_1M_2\cdots M_m$ on $\cH_\cS$.

Finally, we recall the main features of these matrices $M_j$
inherited from those of $e^{i {K_{j}}}$ which shows that they
indeed satisfy properties (1) and (2) of an RRDO. 
The next lemma is an easy generalization of Proposition 2.1 in \cite{bjm}.
\begin{lem}[\cite{bjm}]
\label{contraction} 
With the definition (\ref{defmj}), we have for any $j\in \N^*$, $M_j\psi_\cS=\psi_\cS$. 
Moreover, to any $\phi\in \cH_\cS$ corresponds a unique $A\in \fm_\cS$ such that $\phi=A\psi_\cS$, and $|||\phi|||:=\|A\|_{\cB(\cH_\cS)}$ 
defines a norm on $\cH_\cS$. As an operator on $\cH_\cS$ endowed with the norm $|||\cdot |||$,  
$M_j$ is a contraction for any $j\in\N^*$.
\end{lem}


\subsection{Proof of Theorem \ref{thm4}}

For random repeated interactions, Theorem \ref{thm3} shows that there exists a vector $\theta\in \cH_\cS$
satisfying $\bra \theta, \psi_\cS\ket =1$, such that, with probability $1$, 
\be\label{excm}
\lim_{M\ra\infty}\frac{1}{M}\sum_{m=1}^{M}\Psi_m(\omega)=|\psi_\cS\ket \bra \theta |=\pi,
\ee
where $\pi =|\psi_\cS\ket\bra \theta |$ is a rank one projector in $\cH_\cS$. Therefore, (\ref{almost}) yields  
\bea
&&\hspace{-3cm}\lim_{M\ra\infty}\frac{1}{M}\sum_{m=1}^{M}\bra \psi | (B')^*  B'\alpha_{RI}^m(A_\cS \otimes\one) \psi\ket\nonumber\\
&=&\bra \psi,  (B')^*  B' \tilde U^{+}(N)^* e^{i {K_1}}\cdots 
e^{i{K_{N}}} P  (\pi A_\cS \otimes\one) \psi \ket\nonumber\\
&=& \bra \psi, (B')^*  B' \tilde U^{+}(N)^* e^{i {K_1}}\cdots 
e^{i{K_{N}}}|\psi\ket \  \bra \theta,  A_\cS \psi_\cS\ket_{\cH_\cS}\nonumber\\
&=&\bra \psi | (B')^*  B' \psi\ket\ \bra \theta | A_\cS \psi_\cS\ket_{\cH_\cS}\nonumber\\
&=& \|B' \psi\|^2 \bra\theta | A_\cS \psi_\cS\ket_{\cH_\cS},\nonumber
\eea
where $\|B' \psi\|$ is arbitrary close to $\|\psi\|=1$, for $N$ large enough, see (\ref{psib}). 

In other words, (\ref{excm}) implies the existence of an asymptotic ergodic state $\varrho_E$
on $\fm_\cS$, defined by $\varrho_E(A_\cS):=\bra\theta, A_\cS \psi_\cS\ket_{\cH_\cS}$, $\forall A_\cS\in \fm_\cS$, 
such that for any normal state $\varrho $ on $\fm$,  $\lim_{M\ra\infty}\frac{1}{M}\sum_{m=1}^{M}\varrho(\alpha_{RI}^m(A_\cS \otimes\one))=\varrho_E(A_\cS)$, 
which is the statement of Theorem \ref{thm4}. \hfill \qed


\section{Random Repeated Interaction Spin-Spin system}
\label{sec:spin}

We consider the application of our general results to the concrete situation where both $\cS$ and $\cE$ are two-level ``atoms'', and where the 
interaction induces energy exchange processes. This is a particular case of the
third example in \cite{bjm}. We present the following results
without proofs, a complete analysis of this (and more general) examples will be given in \cite{bjm2}.


\subsection{Description of the model}\label{ssec:spinmodel}

The observable algebra for $\cS$ and for $\cE$ is $\fa_\cS=\fa_\cE=M_2(\C)$. 
Let $E_\cS, E_\cE>0$ be the ``excited'' energy level of $\cS$ and of $\cE$, respectively. Accordingly, the Hamiltonians are given by
$$
h_\cS=
\left[ 
\begin{array}{cc} 0 & 0 \\ 
0 & E_\cS 
\end{array}
\right] 
\mbox{\ \ and \ \ }
h_\cE=
\left[ 
\begin{array}{cc} 
0 & 0 \\ 
0 & E_\cE 
\end{array}
\right].
$$ 
The dynamics are given by $\alpha_\cS^t(A)= \e^{i t h_\cS}A\e^{-i t h_\cS}$  and  $\alpha_\cE^t(A)= \e^{i t h_\cE}A\e^{-i t h_\cE}$. 
We choose the reference state of $\cE$ to be the Gibbs state at inverse temperature $\beta$, i.e., 
\begin{equation}
\varrho_{\beta,\cE}(A)=\frac{\tr(\e^{-\beta h_\cE}A)}{\tr(\e^{-\beta h_\cE})},
\label{mm110}
\end{equation}
and we choose the reference state for $\cS$ to be the tracial state, $\varrho_{0,\cS}(A)=\frac{1}{2}\tr(A)$. The interaction operator is defined by
$$
v :=a_\cS\otimes a_\cE^*+a_\cS^*\otimes a_\cE,
$$
where 
$$
a_{\#}=
\left[
\begin{array}{cc} 
0 & 1 \\ 
0 & 0 
\end{array}
\right]
\mbox{\ \ and\ \ }
a_{\#}^*=
\left[
\begin{array}{cc} 
0 & 0 \\ 
1 & 0 
\end{array}
\right]
$$
are the annihilation and creation operators respectively, of $\#=\cs, \cE$. The Heisenberg dynamics of $\cS$ coupled to one element $\cE$ is given by the $*$-automorphism group $t\mapsto \e^{ith_\lambda} Ae^{-ith_\lambda}$, $A\in{\frak A}_\cS\otimes{\frak A}_\cE$, $h_\lambda = h_\cS+h_\cE+\lambda v$, where $\lambda$ is a coupling parameter.

To find a Hilbert space description of the system, one performs the Gelfand-Naimark-Segal construction of $(\fa_\cS,\omega_\cS)$ and $(\fa_\cE,\omega_\cE)$, see e.g. \cite{BR,bjm}. In this representation, the Hilbert spaces are given by  
$\cH_\cS = \cH_\cE = {\mathbb C}^2\otimes{\mathbb C}^2$,
and the vectors representing $\omega_\cs$ and $\omega_\cE$ are
\begin{equation}
\psi_\cS = \frac{1}{\sqrt 2} \sum_{j=1,2}\varphi_j\otimes\varphi_j \mbox{\ \ and\ \ } 
\psi_\cE=\frac{1}{\sqrt{{\rm Tr}\e^{-\beta h_\cE}}}\left[ \varphi_1\otimes\varphi_1 + \e^{-\beta E_\cE/2}\varphi_2\otimes\varphi_2\right],
\label{mm101}
\end{equation}
respectively, i.e., we have $\varrho_{\beta_{\#},\#}(A)=\scalprod{\psi_{\#}}{(A\otimes\bbbone)\psi_{\#}}$, $\#=\cS,\cE$, $\beta_\cE=\beta$, $\beta_\cS=0$. Here, we have set
$$
\varphi_1 = 
\left[
\begin{array}{c}
1\\
0
\end{array}
\right]
\mbox{\ \ and\ \ }
\varphi_2 = 
\left[
\begin{array}{c}
0\\
1
\end{array}
\right].
$$
We now illustrate our results on this model.


\subsection{Variable or random interaction times}
\label{sssec:time}

The situation where $\tau$ varies is particular in the sense that we can chose $\psi(\omega)\equiv\psi$ so that we don't need to consider ergodic means, see also Remark 8 after Theorem \ref{thm3}. We consider both the deterministic and the random situations. In the deterministic case, we require that the sequence $(\tau_n)_n$ avoids the set $T\N$, where 
$$
T:=\frac{2\pi}{\sqrt{(E_\cS-E_\cE)^2+4\lambda^2}}.
$$
For interaction times in $T\N$, there appear resonance phenomena that require a more refined analysis of the dynamics. This has already been noticed in \cite{bjm}. In the random setting, we only require that the probability of the random variable $\tau(\omega)$ to avoid $T\N$ is non-zero. Let us denote by $\ds \d(\tau,T\N):=\inf_{k\in\N} |\tau-kT|$ the distance between the interaction time $\tau$ and the set $T\N$. 
 
\begin{thm}\label{thm:spintaudeterm} Let $(\tau_n)_n$ be a sequence of interaction times. Assume that there exists $\epsilon>0$ such that $\d(\tau_n,T\Z)>\epsilon$ for all $n$. Then there exists a $\gamma>0$ such that  
$$
\Psi_n=\frac{2}{\tr(\e^{-\tbeta h_\cS})}|\psi_\cS\ket\bra \e^{-\tbeta h_\cS}\otimes\one \ \psi_\cS| +O(\e^{-\gamma n}), 
$$
where $\tbeta=\frac{E_\cE}{E_\cS}\beta$ and where $\psi_\cS$ is given in \fer{mm101}. Consequently, we have for any initial state $\varrho$, and for any $A_\cS\in\fm_\cS$, 
$$
\varrho(\alpha_{RI}^n(A_\cS\otimes\one_\cC))=\varrho_{\tbeta,\cS}(A_\cS)+O(\e^{-\gamma n}),
$$
where $\varrho_{\tbeta,\cS}$ is the Gibbs state \fer{mm110} of the small system for the inverse temperature $\tbeta$.
\end{thm}

Physically, it is reasonable to assume that the interaction time $\tau(\omega)$ takes values in an interval of uncertainty, since during an experiment, the interaction time cannot be controlled exactly, but rather fluctuates around some mean value. We then have the following result.

\begin{thm}
\label{thm:spintaurandom} 
Let $\Omega\ni\omega \to \tau(\omega)\in \R^*_+$ be a random variable. We assume that 
$\p \left(\tau(\omega)\notin T\Z \right)>0$.
Then there exist a set $\Omega_6\subset\Omega^{{\mathbb N}^*}$ and a constant $\gamma>0$ such that ${\mathbb P}(\Omega_6)=1$, and, for all $\po\in\Omega_6$, there exists $C_{\po}$ s.t.
$$
\left\| \Psi_n(\po)-\frac{2}{\tr(\e^{-\tbeta h_\cS})}|\psi_\cS\ket\bra \e^{-\tbeta h_\cS}\otimes\one \ \psi_\cS| \right\| \leq C_{\po}\e^{-\gamma n}.
$$
Consequently, we have for any initial state $\varrho$, and for any $A_\cS\in\fm_\cS$
$$
\varrho(\alpha_{RI}^{n,\po}(A_\cS\otimes\one_\cC))=\varrho_{\tbeta,\cS}(A_\cS)+O_{\po}(\e^{-\gamma n}) \ \ \ a.e.
$$
\end{thm}


\subsection{Random energies of $\cE$}
\label{sssec:atom}

We consider the situation where the energy $\cE$ is a random variable. Set
$$
e_0(\omega):=\frac{\left(E_\cS-E_\cE(\omega)-\sqrt{(E_\cS-E_\cE(\omega))^2+4\lambda^2}\right)^2+4\lambda^2 e^{i\tau \sqrt{(E_\cS-E_\cE(\omega))^2+4\lambda^2}}}{\left(E_\cS-E_\cE(\omega)-\sqrt{(E_\cS-E_\cE(\omega))^2+4\lambda^2}\right)^2+4\lambda^2}.
$$
We show in \cite{bjm2} that this is an eigenvalue of $M(\omega)$, and that $1$ is a degenerate eigenvalue of $M(\omega)$ if and only if 
$$
\tau\frac{\sqrt{(E_\cS-E_\cE(\omega))^2+4\lambda^2}}{2\pi} \in \Z
$$
(in which case $e_0(\omega)=1$). 
In the statement of the next result, we denote the average of a random variable $f(\omega)$ by $\overline f$.

\begin{thm}
\label{thm:spinatom} 
Let $\Omega\ni\omega \to E_\cE(\omega)\in \R^*_+$ be a random variable. We assume that 
\begin{equation}
\p \left( \tau\frac{\sqrt{(E_\cS-E_\cE(\omega))^2+4\lambda^2}}{2\pi} \notin \Z \right)>0.
\label{mm102}
\end{equation}
Then there exists a set $\Omega_7\subset\Omega^{{\mathbb N}^*}$ such that ${\mathbb P}(\Omega_7)=1$, and such that for any $\po\in\Omega_7$,
\be
\lim_{N\to\infty} \frac{1}{N} \sum_{n=1}^N \Psi_n(\po)=|\psi_\cS \ket\bra \rho_E \otimes \one \ \psi_\cS |,\nonumber
\ee
where 
\bea
\rho_E & = & \left[1-(1-\overline{e_0})^{-1}\overline{(1-e_0)(1-2 Z_{\tbeta,\cS}^{-1})}\right] |\varphi_1\rangle\langle\varphi_1| \nonumber\\
 & & \qquad\qquad + \left[1+(1-\overline{e_0})^{-1}\overline{(1-e_0)(1-2 Z_{\tbeta,\cS}^{-1})}\right] |\varphi_2\rangle\langle\varphi_2|. \nonumber
\eea
Consequently, we have for any initial state $\varrho$, for any $A_\cS\in\fm_\cS$, and for any $\po\in\Omega_7$,
$$
\lim_{N\to\infty} \frac{1}{N} \sum_{n=1}^N \varrho(\alpha_{RI}^{n,\po}(A_\cS\otimes\one_\cC))=\varrho_{\tilde{\beta},\cS}(A_\cS),
$$
where $\varrho_{\tilde{\beta},\cS}$ is the Gibbs state \fer{mm110} of the small system for the inverse temperature 
\be\label{averagetemp}
\tilde{\beta}:= \frac{1}{E_\cS} \log \left(2 \left[1-(1-\overline{e_0})^{-1}\overline{(1-e_0)(1-2 Z_{\tbeta,\cS}^{-1})}\right]^{-1}-1 \right).
\ee
\end{thm}

{\it Remarks.\ }
1. The expression (\ref{averagetemp}) for the ``asymptotic temperature'' is of course also valid in the case of Theorem \ref{thm:spintaurandom}. However, in that case, $\tbeta(\omega)\equiv \tbeta$ and $\tilde\beta=\tbeta$.

2. In general, $\tilde{\beta}\neq \overline\tbeta$ which reflects the fact that in Lemma \ref{lem:mbar}, we usually have $\theta\neq\overline\psi$, as mentioned in Remark 2 after the proof of the lemma.

3.  Physically, one can imagine that different kinds of atoms are injected in the cavity, or that from time to time some ``impurity'' occurs. It is therefore reasonable to consider a finite probability space $\Omega$, so that assumption \fer{mm102} is satisfied provided at least one atom has excited energy $E_\cE$ satisfying $\tau\frac{\sqrt{(E_\cS-E_\cE)^2+4\lambda^2}}{2\pi} \notin \Z$.

4. Of course, one can combine Theorems \ref{thm:spintaurandom} and \ref{thm:spinatom} to obtain a result for systems where the interaction times, excitation energies (and other parameters, like temperatures of the ``incoming atoms'') fluctuate.

\bigskip
{\bf Acknowledgements.\ } We thank the ACI Jeunes Chercheurs ``Stochastic modeling of non-equilibrium systems'' for financial support.


\end{document}